\documentclass[12pt]{report} 

\usepackage{amssymb}
\usepackage{latexsym}

\input diagrams \diagramstyle[PS] 

\newtheorem{thm}{Theorem}[chapter]
\newtheorem{algo}[thm]{Algorithm}
\newtheorem{conj}[thm]{Conjecture}
\newtheorem{cor}[thm]{Corollary}
\newtheorem{defn}[thm]{Definition}
\newtheorem{lemma}[thm]{Lemma}

\newcommand{\ba}{\leftarrow}
\newcommand
  {\binom}[2]{{\left(\!\begin{array}{c}{#1}\\{#2}\end{array}\!\right)}}
\newcommand{\C}{{\mathbb C}}
\renewcommand{\d}{{\mathbf d}}
\newcommand{\dsum}{\oplus}
\newcommand{\Dsum}{\bigoplus}
\newcommand{\Dt}{{\Delta_1\ldots\Delta_t}}
\newcommand{\Dtmt}{{\Delta_1^{m_1}\!\!\ldots\Delta_t^{m_t}}}
\newcommand{\End}{{\mathop{\mathrm{End}}\nolimits}}
\newcommand{\g}{{\mathfrak g}}
\newcommand{\ghat}{{\hat{\mathfrak g}}}
\newcommand{\hyt}{{\mathop{\mathrm{ht}}\nolimits}}
\newcommand{\id}{{\mathop{\mathrm{id}}\nolimits}}
\newcommand{\into}{\hookrightarrow}
\newcommand{\iso}{\simeq}
\renewcommand{\l}{\ell}

\newcommand{\oplustag}[1]{\mathbin{\mathop{\oplus}\limits^{#1}}}
\newcommand{\qed}{\hfill $\Box$ \bigskip}
\newcommand{\R}{{\mathbb R}}
\newcommand{\rank}{{\mathop{\mathrm{rank}}\nolimits}}
\newcommand{\Sym}{{\mathop{\mathrm{Sym}}\nolimits}}
\newcommand{\tensor}{\otimes}
\newcommand{\tensorhat}{\mathbin{\widehat{\tensor}}}
\newcommand{\tl}{{T(\l)}}
\newcommand{\tr}{{\mathop{\mathrm{tr}}\nolimits}}

\newcommand{\w}{\omega}
\newcommand{\wl}{{\omega_\l}}
\newcommand{\wml}{{W_m(\l)}}
\newcommand{\wmax}{{\omega_{\mathrm{max}}}}
\newcommand{\Z}{{\mathbb Z}}

\title{Finite Dimensional Representations of Quantum Affine Algebras}
\author{Michael Kleber}

\begin{document}

%

\thispagestyle{empty}
\null\vfil
\begin{center}
    {\bfseries 
Finite Dimensional Representations of Quantum Affine Algebras
    \par}
    \bigskip \medskip
    by \par
    \bigskip \medskip
Michael Kleber
    \par
\vspace{6ex}
    A dissertation submitted in partial satisfaction of the \par
\smallskip
    requirements for the degree of \par
\bigskip \medskip
    Doctor of Philosophy \par
\bigskip \medskip
    in \par
\bigskip \medskip
    Mathematics \par
\bigskip \medskip
    in the \par
\smallskip
    Graduate Division \par
\smallskip
    of the \par
\smallskip
    University of California at Berkeley
\end{center}
{
\vspace{6ex}
  Committee in charge:
\begin{quote}
Professor Nicolai Yu. Reshetikhin, Chair \\
Professor Vera Serganova \\
Professor Deborah Nolan
\end{quote}
  \begin{center}
    Spring 1998
  \end{center}
}
\vfil\null

\pagebreak

%

\newcommand{\abstractsignature}{}
\begin{abstract}

In 1987, Kirillov and Reshetikhin conjectured a formula for how certain 
finite-dimensional representations of the quantum affine algebra 
$U_q(\ghat)$ decomposed into $U_q(\g)$-modules.  Their conjecture was 
built on techniques from mathematical physics and the fact that the 
characters of those particular representations seem to satisfy a certain 
set of polynomial relations, generalizations of the discrete Hirota 
equations.

We present a new interpretation of this formula, involving the 
geometry of weights in the Weyl chamber of $\g$.  This revision has 
the virtue of being computationally easy, especially compared to the 
original form: the Kirillov-Reshetikhin version of the formula is 
computationally intractable for all but the simplest cases.  The 
original version parameterizes the pieces of the decomposition in terms 
of combinatorial objects called ``rigged configurations'' which are 
very hard to enumerate.  We give a bijection between rigged 
configurations and simpler combinatorial objects which can be easily 
generated.

This new version of the formula also adds some structure to the
decomposition: the irreducible $U_q(\g)$-modules are naturally the nodes
of a tree, rooted at the representation containing the original highest
weight vector.  This new tree structure is somewhat consistent among
representations whose highest weights are different multiples of the same
fundamental weight.  We use this coherence of structure to calculate the
asymptotics of the growth of the dimension of these representations as
the multiple of the fundamental weight gets large.

We also explore further the polynomial relations that seem to hold 
among the characters of these representations.  The fact that every 
finite-dimensional representation of the quantum affine algebra is a 
direct sum of representations of the underlying quantized Lie algebra 
is a very strong positivity condition.  We prove that for the 
classical families of Lie algebras, the positivity condition and the 
polynomial relations leave only one choice for the characters of the 
quantum affine algebras --- the ones predicted by the 
Kirillov-Reshetikhin formula.  Therefore to prove the conjecture, it 
would suffice to verify that the characters do indeed satisfy this set 
of relations.

\abstractsignature
\end{abstract}

%

\renewcommand{\thepage}{\roman{page}}
\setcounter{page}{3}

\tableofcontents

\listoffigures
\bigskip\bigskip
\noindent
The pictures of decompositions in Chapter~\ref{chap_decomp} were made
using Paul Taylor's excellent package {\tt diagrams.tex}.  It is
available from any Comprehensive \TeX\ Archive Network (CTAN) site as
\begin{center}
{\tt macros/generic/diagrams/taylor/diagrams.tex}
\end{center}
or directly from its home site,
\begin{center}
{\tt ftp://ftp.dcs.qmw.ac.uk/pub/tex/contrib/pt/diagrams/diagrams.tex}
\end{center}

\pagebreak

\begin{center}
\LARGE\bf Acknowledgements
\end{center}
\bigskip

I am deeply grateful to my advisor, Nicolai Reshetikhin, for helping 
this dissertation come to pass.  Without his guidance and illumination 
I would never have finished; without his patience and insight I might 
never have started.

I am grateful to my wife, Jessica Polito, for support and clear 
thinking.  She has been my constant companion, mathematical and 
otherwise, throughout this adventure and into the next.

Many other members of the Berkeley Math Department, past and present,
contributed to this research over the course of our discussions.  I am
especially glad of having talked with David M\raisebox{.5ex}{c}Kinnon,
David Jones, Ian Grojnowski, Richard Borcherds, and Jim Borger, mostly
for times when they asked me the right questions.  Thanks also to Vera
Serganova for her comments on an earlier draft.

The writing of this dissertation was supported by an Alfred P. Sloan 
Doctoral Dissertation Fellowship.  The research was also partly 
supported by NSF grand DMS 94-01163, and partly conducted while 
visiting the Research Institute for Mathematical Sciences (RIMS), 
Kyoto, Japan, thanks to the generosity of T. Miwa.

\pagebreak
\renewcommand{\thepage}{\arabic{page}}
\setcounter{page}{1}

\pagestyle{headings}

%

\chapter{Introduction}

The theory of finite-dimensional representations of complex simple Lie
algebras is well understood.  Furthermore, if one Lie algebra appears as
a subalgebra of another, due to a corresponding embedding of Dynkin
diagrams, there are well-known branching rules for how representations of
the larger algebra decompose under the action of the smaller one.

Any finite-dimensional complex simple Lie algebra $\g$ is a subalgebra of
its corresponding infinite-dimensional affine Lie algebra $\ghat$.  The
quantized universal enveloping algebra of the affine Lie algebra,
$U_q(\ghat)$, is a Hopf algebra of interest to mathematicians and
mathematical physicists, introduced simultaneously by Drinfeld and Jimbo
around 1985.  Finite-dimensional representation of $U_q(\ghat)$ are not
well understood, and even their structure when viewed as representations
of the Hopf subalgebra $U_q(\g)$ is not generally known.

In 1987, Kirillov and Reshetikhin conjectured a formula for how some
finite-dimensional representations of $U_q(\ghat)$ decomposed into
$U_q(\g)$-modules.  They looked only at representations whose highest
weight is a multiple of a fundamental weight.  Their conjecture was built
on techniques from mathematical physics and the fact that the characters
of those particular representations seem to satisfy a certain set of
polynomial relations, generalizations of the discrete Hirota equations.

In Chapter~\ref{chap_decomp}, we give a new interpretation of this 
formula, involving the geometry of weights in the Weyl chamber of 
$\g$.  This revision has the virtue of being computationally easy, 
especially compared to the original form: the Kirillov-Reshetikhin 
version of the formula is computationally intractable for all but the 
simplest cases.  The original version parameterized the pieces of the 
decomposition in terms of combinatorial objects called ``rigged 
configurations'' which are very hard to enumerate.  We give a 
bijection between rigged configurations and simpler combinatorial 
objects which can be easily generated.

This new version of the formula also adds some structure to the
decomposition: the irreducible $U_q(\g)$-modules are naturally the nodes
of a tree, rooted at the representation containing the original highest
weight vector.  This new tree structure is somewhat consistent among
representations whose highest weights are different multiples of the same
fundamental weight.  We use this coherence of structure to calculate the
asymptotics of the growth of the dimension of these representations as
the multiple of the fundamental weight gets large.

In Chapter~\ref{chap_poly}, we explore further the polynomial relations
that seem to hold among the characters of these representations.  The
fact that every finite-dimensional representation of the quantum affine
algebra is a direct sum of representations of the underlying quantized Lie
algebra is a very strong positivity condition.  We prove that for the 
classical families of Lie algebras, the positivity condition and the 
polynomial relations leave only one choice for the characters of the 
quantum affine algebras --- the ones predicted by the 
Kirillov-Reshetikhin formula.  Therefore to prove the conjecture, it 
would suffice to verify that the characters do indeed satisfy this set 
of relations.

Chapter~\ref{chap_questions} lists some natural questions for further
research.  Mostly, they ask for generalizations of the notions mentioned
above to other contexts, some straightforward and some completely
open-ended.

\chapter[Overview]{Overview of Yangians and Quantum Affine Algebras}
\label{chap_overview}

\section{The Algebras}
\label{sec_algebras}

To any finite-dimensional simple Lie algebra $\g$, we can associate 
two closely-related Hopf algebras: its Yangian $Y(\g)$ and its quantum 
affine algebra $U_q(\ghat)$.

\subsection{The Yangian}

The Yangian was introduced by Drinfeld in \cite{Dr} as part of the study
of solutions to the Quantum Yang-Baxter Equation (this connection is
discussed in section~\ref{sec_QYBE}).  A second definition of Yangians in
terms of generators and relations, with an easier description of the
action on highest-weight modules, was given in \cite{Dr:new}, and this is
the one we give here.

Fix a complex simple Lie algebra $\g$ with simple roots $\alpha_1, 
\ldots, \alpha_r$, $r=\rank(\g)$ with respect to some chosen Cartan 
subalgebra.  Let $C=(c_{ij})$ denote the Cartan matrix of $\g$, and 
let $b_{ij}=(\alpha_i,\alpha_j)/2$ be the symmetrized version.

\begin{defn}
The Yangian $Y(\g)$ is an associative algebra with generators 
$\kappa_{ik}$, $\xi^+_{ik}$, $\xi^-_{ik}$, where $i=1,\ldots,r$ and 
$k=0,1,2,\ldots$, and relations
\begin{eqnarray} &
 [\kappa_{ik},\kappa_{jl}]=0, \,\,
 [\kappa_{i0},\xi^\pm_{jl}]=\pm 2b_{ij}\xi^\pm_{jl}, \,\,
 [\xi^+_{ik},\xi^-_{jl}]=\delta_{ij}\kappa_{i,k+l},
\\ & \label{rel1}
 [\kappa_{i,k+1},\xi^\pm_{jl}] - [\kappa_{ik},\xi^\pm_{j,l+1}] =
   \pm b_{ij}(\kappa_{ik}\xi^\pm_{jl} + \xi^\pm_{jl}\kappa_{ik}),
\\ &\label{rel2}
 [\xi^\pm_{i,k+1},\xi^\pm_{jl}] - [\xi^\pm_{ik},\xi^\pm_{j,l+1}] =
   \pm b_{ij}(\xi^\pm_{ik}\xi^\pm_{jl} + \xi^\pm_{jl}\xi^\pm_{ik}),
\\ &
 \mbox{for $i\neq j$, $n=1-a_{ij}$, }
   \Sym[\xi^\pm_{i,k_1},[\xi^\pm_{i,k_2},\cdots[\xi^\pm_{i,k_n},
     \xi^\pm_{jl}]\cdots]]=0
\end{eqnarray}
where $\Sym$ is the sum over all permutations of $k_1,\ldots,k_n$.
\end{defn}

The action of $Y(\g)$ on finite-dimensional representations is similar 
to the situation for $\g$ itself: in any finite-dimensional module 
$V$, there is a nonzero ``highest weight'' vector $v$, unique up to 
multiplication by scalars, which is sent to $0$ by all the 
$\xi^+_{ik}$ and which is an eigenvector for all $\kappa_{ik}$.  All 
of $V$ is generated by the action of the $\xi^-_{ik}$ on $v$, and the 
analog of the Poincar\'e-Birkhoff-Witt theorem holds, allowing us to 
pick a total order on the generators such that ordered words form a 
linear basis for $Y(\g)$.

In the case of representations of $\g$, an irreducible highest weight module 
is finite-dimensional if the eigenvalues of the action on the highest 
weight vector are nonnegative integers.  There is an analog in 
$Y(\g)$.  Suppose $\kappa_{ik}v=d_{ik}v$, for some $d_{ik}\in\C$.  
Then the module $V$ is finite-dimensional if and only if
$$
1+\sum_{k=0}^\infty d_{ik}u^{-k-1} = \frac{P_i(u+b_{ii})}{P_i(u)}
$$
where $P_i(u)$ is a polynomial in $u$ and the left-hand side is the 
rational function's Taylor series at infinity.  The set of polynomials 
$P_1(u),\ldots,P_r(u)$ are only defined up to a scalar, so we choose 
them to be monic.  They are called the Drinfeld polynomials of $V$; 
$r$-tuples of monic polynomials are in bijection with irreducible 
finite-dimensional $Y(\g)$-modules in this way.  However, given a set 
of Drinfeld polynomials, there is no known way to calculate a 
character or even the dimension of the associated representation.

Note that there is a copy of $\g$ actually embedded in $Y(\g)$, as the 
subalgebra generated by the $\kappa_{i0}$ and $\xi^\pm_{i0}$.  The 
highest weight vector of a $Y(\g)$-module is therefore a highest 
weight vector for an action of $\g$ on the same space.  The highest 
weights of the resulting $\g$ action on a $Y(\g)$-module are exactly 
the degrees of the Drinfeld polynomials $P_i(u)$.

If we multiply the right-hand side of the relations~(\ref{rel1})
and~(\ref{rel2}) by $h$, we get defining relations for another Hopf
algebra, $Y_h(\g)$.  It is a deformation of the loop algebra of
polynomial maps $\C^\times\to\g$ with the pointwise bracket
(see~\cite{ChPbook} for an introduction to deformation and quantization
of Hopf algebras).  In the classical limit $h\to 0$, the generators
$\kappa_{ik}$, $\xi^+_{ik}$, $\xi^-_{ik}$ of $Y_h(\g)$ are sent to the
polynomial loops $H_iu^k$, $X_iu^k$ and $Y_iu^k$ in an indeterminate $u$.
For all values of $h$ other than $h=0$, though, $Y_h(\g)$ is 
isomorphic; this is why we can choose to specialize to $h=1$ and just 
work with $Y(\g)$, as above.

\subsection{Quantum Affine Algebras}

We now turn our attention to quantum affine algebras.  These were 
introduced simultaneously by Drinfeld and Jimbo, also as part of the
pursuit of solutions of the Quantum Yang-Baxter Equation.

The quantum affine algebra can be realized in several different ways.  
First, we let $\ghat$ denote the (untwisted) affine Kac-Moody algebra 
associated to the extended Dynkin diagram of $\g$ (with the added node 
numbered $0$ and corresponding root $\alpha_0$).  In \cite{Dr}, 
Drinfeld showed that the universal enveloping algebra $U(\ghat)$, and 
indeed the universal enveloping algebra of any symmetrizable Kac-Moody 
algebra, can be quantized to give $U_h(\ghat)$.

The algebras $U_h(\ghat)$ are better-understood than for arbitrary 
Kac-Moody algebras because they have a second realization in terms of 
central extension of the loop algebra $\g[u,u^{-1}]$ of Laurent 
polynomial maps $\C^\times\to\g$.  In \cite{Dr:new} Drinfeld provided 
a new definition of $U_h(\g)$, which we copy here, whose generators 
make the loop algebra structure visible: one can think of $\kappa_{ik}$, 
$\xi^+_{ik}$ and $\xi^-_{ik}$ as $H_iu^k$, $X_iu^k$ and $Y_iu^k$.

\begin{defn}
The quantum affine algebra $U_h(\ghat)$ is an $h$-adically complete 
associative algebra over $\C[[h]]$ with generators $\kappa_{ik}$, 
$\xi^+_{ik}$, $\xi^-_{ik}$ and the central element $c$, where 
$i=0,\ldots,r$ and $k\in\Z$, and relations
\begin{eqnarray*} &
 [c,\kappa_{ik}] = [c,\xi^\pm_{ik}] = 0,
\\ &
 [\kappa_{ik},\kappa_{jl}] = 4\delta_{k,-l}k^{-1}h^{-2}
   \sinh(khb_{ij}) \sinh(khc/2),
\\ &
 [\kappa_{ik},\xi^\pm_{jl}] = \pm 2(kh)^{-1}
   \sinh(khb_{ij}) \exp(\mp|k|hc/4) \xi^\pm_{j,k+l},
\\ &
 \xi^\pm_{i,k+1}\xi^\pm_{jl} - e^{\pm hb_{ij}} \xi^\pm_{jl}\xi^\pm_{i,k+1} =
   e^{\pm hb_{ij}}\xi^\pm_{ik}\xi^\pm_{j,l+1} - \xi^\pm_{j,l+1}\xi^\pm_{ik},
\\ &
 [\xi^+_{ik},\xi^-_{jl}] = \delta_{ij}h^{-1}
  \{\psi_{i,k+l}\exp(hc(k-l)/4)-\phi_{i,k+l}\exp(hc(l-k)/4)\}
\\ &
 \mbox{for $i\neq j$, $n=1-a_{ij}$, } \displaystyle
   \Sym\sum_{r=0}^n(-1)^r C_n^r(hb_{ii}/2)
    \xi^\pm_{ik_1} \cdots \xi^\pm_{ik_r} \xi^\pm_{jl}
    \xi^\pm_{ik_{r+1}} \cdots \xi^\pm_{ik_n} = 0
\end{eqnarray*}
where $\Sym$ is the sum over all permutations of $k_1,\ldots,k_n$,
$$
C_n^r(\alpha) = \frac
{\sinh n\alpha \cdot \sinh(n-1)\alpha \cdot\ldots\cdot \sinh(n-r+1)\alpha}
{\sinh\alpha \cdot \sinh2\alpha \cdot\ldots\cdot \sinh r\alpha}
$$
and the $\phi_{ip}$ and $\psi_{ip}$ are determined by the relations
\begin{eqnarray*}
\sum_p \phi_{ip}u^{-p} &=& \exp\left\{ -h\left(
 \frac{\kappa_{i0}}{2}+\sum_{p<0} \kappa_{ip}u^{-p} \right) \right\}, \\
\sum_p \psi_{ip}u^{-p} &=& \exp\left\{ h\left(
 \frac{\kappa_{i0}}{2}+\sum_{p>0} \kappa_{ip}u^{-p} \right) \right\}.
\end{eqnarray*}
\end{defn}

It is frequently more convenient to talk about the deformation 
$U_q(\ghat)$ instead of $U_h(\ghat)$, where $q=e^h$.  Technically, 
$U_h(\ghat)$ is defined over $\C[[h]]$, which forces us to worry 
about $h$-adic completions of algebras and requires careful thinking 
to specialize $h$ to any specific nonzero value.  By looking at 
$U_q(\ghat)$ instead, we can deal with an algebra defined over 
$\C[q,q^{-1}]$ with $q$ a formal variable, and can specialize $q$ to 
any nonzero complex number easily.  If $q$ is a root of unity we get 
different behavior, corresponding to the change in $U_h(\ghat)$ if 
$h$ were nilpotent, but for generic $q$ everything we want about 
$U_h(\ghat)$ is preserved.

\subsection{Equivalence of Decomposition}

Finite-dimensional irreducible representations of $Y(\g)$ and of
$U_q(\ghat)$ are closely related.  In each case they are indexed by Drinfeld
polynomials $P_1,\ldots,P_r$.  As mentioned above, the degrees of the
Drinfeld polynomials in the Yangian case give the highest weight of the
$\g$ action on the module.  Similarly, in the $U_q(\ghat)$ context, the
degrees give the highest weight of the action of $U_q(\g)$, which sits as
a subalgebra inside of $U_q(\ghat)$ based on the inclusion of Dynkin
diagrams.  In fact, it seems these two situations are identical:

\begin{conj}
Let $P_1,\ldots,P_r$ be monic polynomials.  Decompose the $Y(\g)$ module
with those Drinfeld polynomials into $\g$-modules as $\bigoplus
V_\lambda^{\oplus n_\lambda}$.  Likewise, decompose the $U_q(\ghat)$
module with the same Drinfeld polynomials into $U_q(\g)$-modules as
$\bigoplus V_\lambda^{\oplus m_\lambda}$.  Then for each highest weight
$\lambda$, the multiplicities are the same: $n_\lambda = m_\lambda$.
\end{conj}

This appears to be a fact that everyone believes, but no one has
provided a proof.  Statements made in one of these two contexts have
been happily transferred to the other in the literature with no comment.
We regretfully continue to sweep this omission under the rug.

\section{The Yang-Baxter Equation and the Bethe Ansatz}
\label{sec_QYBE}

Yangians and quantum affine algebras originally arose in the study of
mathematical physics.  Most of what are now the axioms of a Hopf algebra
started as as hoc tools for finding solutions to the Quantum Yang-Baxter
equation.  The ties between the two fields were the main subject of 
the paper \cite{Dr} in which Yangians were originally defined.

\subsection{Yang-Baxter Equations}

The Quantum Yang-Baxter Equation (QYBE) is the
following requirement on a matrix $R\in\End(V\tensor V)$:
\begin{equation}
\label{QYBE}
R_{12} R_{13} R_{23} = R_{23} R_{13} R_{12}
\end{equation}
The equation holds in $V\tensor V\tensor V$, where $R_{ij}$ indicates 
that $R$ is acting on the $i$th and $j$th components in the tensor 
product.  Such an $R$ is a constant solution of the more general 
Quantum Yang-Baxter Equation with spectral parameters,
\begin{equation}
\label{QYBEspec}
R_{12}(u-v) R_{13}(u-w) R_{23}(v-w) = R_{23}(v-w) R_{13}(u-w) R_{12}(u-v)
\end{equation}
a functional equation for matrix-valued functions $R:\C\to\End(V\tensor V)$

In any Hopf algebra $A$, there are two possible comultiplications 
$\Delta:A\to A\tensor A$ and its opposite $\Delta^{\mathrm 
op}=\sigma\circ\Delta$, where $\sigma$ acts on $A\tensor A$ by 
switching the factors.  In a cocommutative Hopf algebra the two 
comultiplications are equal.  We say $A$ is {\em almost cocommutative} 
if they are instead conjugate; that is, if there exists an invertible 
element $R$ such that $\Delta^{\mathrm op}=R\Delta R^{-1}$.  In many 
cases there is no such element $R$ in $A\tensor A$ but there is in 
some completion $A\tensorhat A$, which still suits our needs as 
long as conjugation by $R$ stabilizes $A\tensor A$ inside 
$A\tensorhat A$.

We further say that $A$ is {\em quasitriangular} if 
$(\Delta\tensor\id)R = R_{13}R_{23}$ and $(\id\tensor\Delta)R = 
R_{13}R_{12}$, and we call $R$ the {\em universal $R$-matrix} of the 
Hopf algebra.  One can easily check that the universal $R$ matrix is 
automatically a solution to the QYBE.

In the case of Yangians, we get solutions to the QYBE with spectral 
parameters.  While there is no $R$-matrix in $Y(\g)\tensor Y(\g)$ 
itself, there is an element $R(u)$ in the completion $(Y(\g)\tensor 
Y(\g))[[u^{-1}]]$ which has the form
$$
R(u) = 1\tensor1 +\frac{t}{u} + \sum_{n=2}^{\infty} \frac{R_n}{u^n}
$$
which intertwines the comultiplications.  This Taylor series is 
ill-suited for solving the QYBE with spectral parameters directly, 
since equation~(\ref{QYBEspec}) would require multiplying expansions in 
different indeterminates.  Fortunately, one can show that if we let 
$R(u)$ act on any finite-dimensional representation, $R(u) = f(u) 
R_{\mathrm{rat}}(u)$, where $R_{\mathrm{rat}}(u)$ is a rational 
function of $u$ and $f(u)$ is meromorphic away from a countable set of 
points in $\C$.

Therefore finite-dimensional representations of $Y(\g)$ give rise to 
so-called rational solutions of the QYBE with spectral parameters.

\subsection{Transfer matrices and the Bethe Ansatz}

In the study of integrable lattice models, a central role in
understanding the behavior of the system is played by a linear operator
$t(u)$ acting on the space ${\mathcal H} = V_1\tensor
V_2\tensor\cdots\tensor V_n$.  This operator is called the {\em
row-to-row transfer matrix}, and is defined as
$$
t(u) = \tr_0 \, R_{01}(u-w_1) R_{02}(u-w_2) \ldots R_{0n}(u-w_n)
$$
where the factors $R_{0i}$ act in $V_0\tensor{\mathcal H}$ on $V_0$ and
$V_i$, and the resulting operator acts in $\mathcal H$.  The system is
called {\em integrable} if $R(u)$ is a nontrivial solution of the Quantum
Yang-Baxter Equation with spectral parameters~(\ref{QYBEspec}).  Using
the QYBE, one can easily verify that $[t(u),t(v)]=0$.  In this case the
transfer matrix is a generating function for the commuting quantum
Hamiltonians of the associated system.  The spectrum of this commuting
family determines the major characteristics of the system; for an
introduction to statistical mechanics and quantum integrable systems,
see~\cite{ChPbook}.

The Bethe Ansatz is the main technique for calculating the eigenvalues 
of the transfer matrices, pioneered by H. Bethe in the 1930s.  The 
eigenvalues of the transfer matrices are given as the solutions to a 
set of algebraic equations, the Bethe equations.  The method only 
locates eigenvalues corresponding to Bethe vectors, eigenvectors which 
satisfy a certain technical condition.  However, there is evidence 
that finding the Bethe vectors should suffice.  This was the grounds 
for the conjectures we mention below.

\subsection{Decomposition of the Tensor Product}

Now consider the case ${\mathcal H}=V_1\tensor V_2\tensor\cdots\tensor
V_n$ where the $V_i$ are all finite-dimensional $Y(\g)$-modules, and we
use the Yangian $R$-matrix to define the transfer matrix $t(u)$.  The
heart of the connection between the QYBE and representation theory is as
follows:

\begin{thm}
View ${\mathcal H}$ as a $\g$-module, by letting the copy of $\g$
embedded in $Y(\g)$ act diagonally on the tensor product.  Then the
transfer matrix $t(u)$ commutes with the $\g$ action.
\end{thm}

Therefore every eigenspace of the transfer matrices is a sum of
$\g$-modules.

\begin{conj}
The spectrum of the transfer matrix is simple with respect to the
$\g$-action.
\end{conj}

This is the best possible scenario.  In this case the action of the
transfer matrix would completely decompose the tensor product into
irreducibles, and we would have a bijection between eigenvalues of $t(u)$
and highest weight vectors of ${\mathcal H}$.  In particular,
Conjecture~\ref{KRconj} in the next section is precisely the statement
that the multiplicity of an irreducible $\g$-module in the tensor product
is just the number of eigenvectors that are highest weight vectors with
the correct highest weight.

It has been proved that the spectrum is simple in some cases.  When 
$V_i\iso\C^n$, the $R$-matrix comes from $Y(\mathfrak{sl}_n)$, and the 
spectral parameters $w_1,\ldots,w_n$ are generic, it was proved 
in~\cite{Ki} that the Bethe vectors lead to a representation of the 
correct dimension.  A bijection between the Bethe vectors and the 
irreducible pieces of the decomposition was completed in~\cite{KKR} 
and~\cite{KR1}.  There is considerable computational evidence, 
including the decompositions in section~\ref{sec_table} here, that the 
conjecture is true in general.

\section{A Result of Kirillov and Reshetikhin}

In \cite{KR}, Kirillov and Reshetikhin used the correspondence between
irreducible $\g$-modules in a $Y(\g)$-module and solutions to the Bethe
equations, along with techniques from mathematical physics, to conjecture
a formula for the decomposition of certain representations of $Y(\g)$.
Since it is sometimes unclear which statements are conjectures and 
which are theorems, we give a precise account of the results from 
that paper in this section.

First, we restrict our attention to representations of Yangians which are
tensor products of $Y(\g)$-modules whose highest weights (when viewed as
$\g$-modules) are multiples of a fundamental weight.  Write
$\alpha_1,\ldots,\alpha_r$ for the fundamental roots and
$\omega_1,\ldots,\omega_r$ for the fundamental weights; $\wml$ is a
$Y(\g)$-module with highest weight $m\wl$, for some $m\in\Z_+$ and
$1\leq\l\leq r$ (see section~\ref{sec_decomp_intro} for precise 
definitions).  We want to decompose
\begin{equation}
\label{over-decomp}
\bigotimes_{a=1}^N (W_{m_a}(\l_a) |_{\g})
\iso \Dsum_{\lambda} V_\lambda^{\dsum n_\lambda}
\end{equation}
The sum runs over all weights $\lambda$ less than $\sum m_a\w_{\l_a}$,
the highest weight of the tensor product.  The nonnegative integer
$n_\lambda$ is the multiplicity with which the irreducible $\g$-module
$V_\lambda$ with highest weight $\lambda$ occurs in the decomposition.

The main result of \cite{KR} is the following:

\begin{conj}[Kirillov-Reshetikhin]
\label{KRconj}
Write $\lambda =  \sum m_a\w_{\l_a} - \sum n_i \alpha_i$.  Then
$$
n_\lambda = 
\sum_{\mbox{partitions}} \;\; \prod_{n\geq1} \;\; \prod_{k=1}^r
\binom{P^{(k)}_n(\nu) + \nu^{(k)}_n}{\nu^{(k)}_n}
$$
The sum is taken over all ways of choosing partitions
$\nu^{(1)},\ldots,\nu^{(r)}$ such that $\nu^{(i)}$ is a partition of
$n_i$ which has $\nu^{(i)}_n$ parts of size $n$ (so $n_i =
\sum_{n\geq1} n \nu^{(i)}_n$).  The function $P$ is defined by
\begin{eqnarray*}
P^{(k)}_n(\nu) &=& \sum_{a=1}^N \min(n,m_a)\delta_{k,\l_a}
  - 2 \sum_{h\geq 1} \min(n,h)\nu^{(k)}_{h} + \\
&&\hspace{1cm} +
  \sum_{j\neq k}^r \sum_{h\geq 1} \min(-c_{k,j}n,-c_{j,k}h)\nu^{(j)}_{h}
\nonumber
\end{eqnarray*}
where $C=(c_{i,j})$ is the Cartan matrix of $\g$, and ${a\choose b}=0$
whenever $a<b$.
\end{conj}
Earlier papers~\cite{KKR} and~\cite{KR1} gave a purely combinatorial
proof of this formula in the case $\g=\mathfrak{sl}_n$, where the sets of
partitions which lead to nonzero binomial coefficients are called {\em
rigged configurations}.  The formula is inspired by counting solutions to
the Bethe equations.  These solutions form ``strings'' and ``holes'': the
numbers $\nu_n^{(k)}$ are the number of color $k$ strings of length $n$,
and the formula for $P^{(k)}_n(\nu)$ counts the corresponding number of
holes.

While this formula is meant to apply to all complex simple Lie algebras, 
the remainder of the paper restricts its attention to the classical 
cases.

First, it is noted that the Yangian is known to act on the following 
spaces:
$$
\begin{array}{rlcll}
A_n:
 & W_1(\l) &=& V(\w_\l) & 1\leq\l\leq n \\
B_n:
 & W_1(\l) &=& V(\w_\l) \oplus V(\w_{\l-2}) \oplus V(\w_{\l-4})
               \oplus\cdots & 1\leq\l\leq n-1 \\
 & W_1(n)  &=& V(\w_n) \\
C_n:
 & W_1(\l) &=& V(\w_\l) & 1\leq\l\leq n \\
D_n:
 & W_1(\l) &=& V(\w_\l) \oplus V(\w_{\l-2}) \oplus V(\w_{\l-4})
               \oplus\cdots & 1\leq\l\leq n-2 \\
 & W_1(\l) &=& V(\w_\l) & \l=n-1,n
\end{array}
$$
While nontrivial to check, these decompositions are indeed the same as 
the ones predicted by equation~(\ref{over-decomp}) and 
Conjecture~\ref{KRconj} in the special case that $N=1$ and $m_1=1$.

Having shown the formula is true for the obvious base cases, one might 
hope to complete a proof of the conjecture by induction.  In the case 
$\g=\mathfrak{sl}_n$, the characters of $\wml$ are known to satisfy a 
certain set of quadratic recurrence relations.  In an earlier 
paper~(\cite{Ki}), Kirillov showed that the characters predicted for 
$\g=\mathfrak{sl}_n$ by Conjecture~\ref{KRconj} also satisfy those 
same recurrence relations, using the Littlewood-Richardson rule.  
Since the base cases just mentioned are a complete set of ``initial 
data'' for the recurrence, this completed the proof of the conjecture 
in the $A_n$ case.

The remainder of~\cite{KR} gives a generalization of half of this 
proof.  The authors write down a set of recurrence relations 
generalizing those known for $A_n$ to the $B_n$, $C_n$ and $D_n$ cases 
(see section~\ref{sec_relations} for these and a version which covers 
the exceptional Lie algebras as well).  Then, although the gruesome 
combinatorial details do not appear in the paper, they verify that the 
characters predicted by the conjecture obey these recurrence 
relations.

To prove the conjectural formulas, it only remains to show that the 
actual characters of the $Y(\g)$ modules $\wml$ indeed satisfy 
these generalized recurrence relations.  Unfortunately, no proof of 
this fact is currently known.

\subsection{Practical Questions}

The formula for $n_\lambda$ given in Conjecture~\ref{KRconj} has one
major flaw: practical computation with this formula is impossible for all
but the simplest examples.

Recall that the formula is a summation over all partitions of a product
of binomial coefficients.  The binomial coefficient ${a\choose b}$ is
defined to be zero whenever $a<b$, which happens any time $P_n^{(k)}(\nu)$
is negative.  So the nonzero terms in the summation correspond to choices
of partitions $\nu^{(1)},\ldots,\nu^{(r)}$ which have the property that
$P_n^{(k)}$ is nonnegative for all $1\leq k\leq r$ and for all
$n=1,2,3,\ldots$.

As the integers $n_1,\ldots,n_r$ in $\lambda = \sum m_a\w_{\l_a} - \sum
n_i \alpha_i$ get larger, the total number of partitions grows much more
quickly than the number which yield nonzero terms in the sum.  Even in
the case of the fundamental representations $W_1(\l)$, where many of the
decompositions were known using other techniques, this problem made it
impossible to verify that the conjecture gave the correct results.

As a practical example, suppose one wanted to calculate the multiplicity
of the trivial representation in $W_1(4)$ for $E_8$ (where $\w_4$
corresponds to the trivalent node of the Dynkin diagram).  The integers
$n_1,\ldots,n_8$ are the $\alpha$-coordinates of $\w_4$, $(10, 15, 20,
30, 24, 18, 12, 6)$, and the number of possible choices for
$\nu^{(1)},\ldots,\nu^{(r)}$ is the product of their partition numbers,
$13,339,892,309,691,024,000$.

One goal of Chapter~\ref{chap_decomp} is to overcome this difficulty.  
Using the methods there, we find that of those $13$ quintillion 
choices, exactly six give nonzero summands, and the total multiplicity 
is ten.

\begin{figure}
\caption{Numbering of nodes on Dynkin diagrams}
\label{table_dd}
$$
\setlength{\unitlength}{1cm}
\thicklines
\newcommand{\ci}{\circle*{.13}}
\renewcommand{\mp}{\multiput}
\newcommand{\num}[1]{{\raisebox{-.5\unitlength}{\makebox(0,0)[b]{${#1}$}}}}
\begin{array}{lll}
\\
A_n && 
\raisebox{1ex}{\begin{picture}(5,0)
\mp(0,0)(1,0){3}{\ci}
\put(0,0){\num{1}}\put(1,0){\num{2}}\put(2,0){\num{3}}
\put(3,0){\makebox(0,0){\ldots}}
\mp(5,0)(-1,0){2}{\ci}
\put(4,0){\num{n-1}}\put(5,0){\num{n}}
\put(0,0){\line(1,0){2.5}}
\put(5,0){\line(-1,0){1.5}}
\end{picture}}
\\ \\
B_n && 
\raisebox{1ex}{\begin{picture}(5,0)
\mp(0,0)(1,0){3}{\ci}
\put(0,0){\num{1}}\put(1,0){\num{2}}\put(2,0){\num{3}}
\put(3,0){\makebox(0,0){\ldots}}
\mp(5,0)(-1,0){2}{\ci}
\put(4,0){\num{n-1}}\put(5,0){\num{n}}
\put(0,0){\line(1,0){2.5}}
\put(4,0){\line(-1,0){.5}}
\put(4,.06){\line(1,0){1}}
\put(4,-.06){\line(1,0){1}}
\put(4.5,0){\makebox(0,0){\Large$>$}}
\end{picture}}
\\ \\
C_n && 
\raisebox{1ex}{\begin{picture}(5,0)
\mp(0,0)(1,0){3}{\ci}
\put(0,0){\num{1}}\put(1,0){\num{2}}\put(2,0){\num{3}}
\put(3,0){\makebox(0,0){\ldots}}
\mp(5,0)(-1,0){2}{\ci}
\put(4,0){\num{n-1}}\put(5,0){\num{n}}
\put(0,0){\line(1,0){2.5}}
\put(4,0){\line(-1,0){.5}}
\put(4,.06){\line(1,0){1}}
\put(4,-.06){\line(1,0){1}}
\put(4.5,0){\makebox(0,0){\Large$<$}}
\end{picture}}
\\ \\
D_n && 
\raisebox{1ex}{\begin{picture}(5,.8)
\mp(0,0)(1,0){3}{\ci}
\put(0,0){\num{1}}\put(1,0){\num{2}}\put(2,0){\num{3}}
\put(3,0){\makebox(0,0){\ldots}}
\put(4,0){\ci}
\put(4,0){\num{n-2\,\,\,}}
\put(0,0){\line(1,0){2.5}}
\put(4,0){\line(-1,0){.5}}
\put(4,0){\line(2,1){1}}
\put(4,0){\line(2,-1){1}}
\put(5,.5){\ci}\put(5,-.5){\ci}
\put(5,.5){\makebox(0,0)[l]{\hspace{.3\unitlength}$n-1$}}
\put(5,-.5){\makebox(0,0)[l]{\hspace{.3\unitlength}$n$}}
\end{picture}}
\\ \\
E_{6,7,8} && 
\raisebox{1ex}{\begin{picture}(5,1.2)
\mp(0,0)(1,0){7}{\ci}\put(2,1){\ci}
\put(0,0){\num{1}}\put(1,0){\num{3}}\put(2,0){\num{4}}\put(3,0){\num{5}}
\put(4,0){\num{6}}\put(5,0){\num{7}}\put(6,0){\num{8}}
\put(2,.9){\makebox(0,0)[r]{$2$\hspace{.2\unitlength}}}
\put(2,0){\line(0,1){1}}
\put(0,0){\line(1,0){4}}
\put(5,0){\line(-1,0){.7}}
\put(6,0){\line(-1,0){.7}}
\end{picture}}
\\ \\
F_4 && 
\raisebox{1ex}{\begin{picture}(5,0)
\mp(1,0)(1,0){4}{\ci}
\put(1,0){\num{1}}\put(2,0){\num{2}}\put(3,0){\num{3}}\put(4,0){\num{4}}
\put(1,0){\line(1,0){1}}
\put(3,0){\line(1,0){1}}
\put(2,.06){\line(1,0){1}}
\put(2,-.06){\line(1,0){1}}
\put(2.5,0){\makebox(0,0){\Large$>$}}
\end{picture}}
\\ \\
G_2 && 
\raisebox{1ex}{\begin{picture}(5,0)
\mp(2,0)(1,0){2}{\ci}
\put(2,0){\num{1}}\put(3,0){\num{2}}
\put(2,0){\line(1,0){1}}
\put(2,.06){\line(1,0){1}}
\put(2,-.06){\line(1,0){1}}
\put(2.5,0){\makebox(0,0){\Large$<$}}
\end{picture}}
\\
\end{array}
$$
\end{figure}
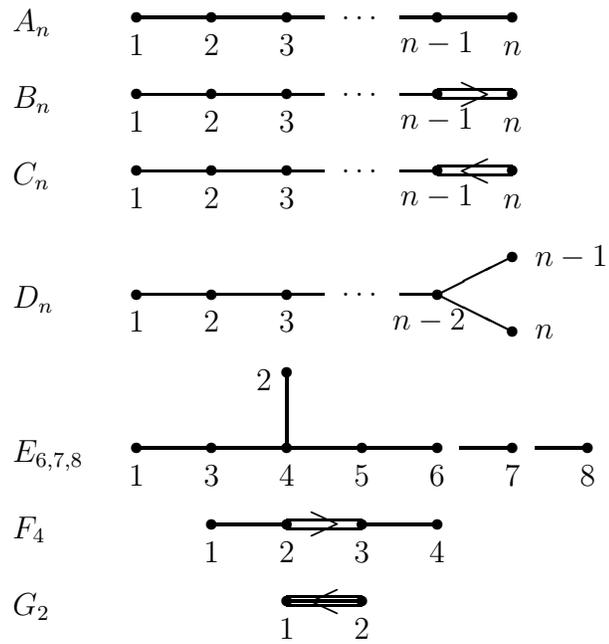

\chapter{Combinatorics of Decomposition}
\label{chap_decomp}

In this chapter we investigate Kirillov and Reshetikhin's conjectured 
formula~\cite{KR} for decomposing certain representations of Yangians 
into irreducible $\g$-modules.  We develop a practical way to compute 
this decomposition, and find some new structure to these modules.  As 
a special case, we can decompose into irreducibles the tensor product 
of an arbitrary number of representations of $\mathfrak{sl}_n$ whose 
associated Young diagrams are rectangles, in a way which is symmetric 
in all the factors.

A preliminary version of this chapter was published in~\cite{kleber}.
The material has been reorganized and some changes have been made
throughout.  In particular, that version did not go into detail about
tensor products of representations.  Section~\ref{sec_rectangles} is new.

\section{Introduction}
\label{sec_decomp_intro}

Let $\g$ be a complex semisimple Lie algebra of rank $r$ and $Y(\g)$ 
its Yangian, as in section~\ref{sec_algebras}.  Write 
$\alpha_1,\ldots,\alpha_r$ for the fundamental roots and 
$\omega_1,\ldots,\omega_r$ for the fundamental weights of $\g$.
We normalize the Killing form so that the long roots
have length $2$.

\begin{defn}
\label{def_wml}
For $\l=1,2,\ldots,r$ and $m=0,1,2,\ldots$, let $\wml$ denote the
irreducible $Y(\g)$-module with Drinfeld polynomials
\begin{eqnarray*}
P_\l(u) &=& \prod_{i=1}^{m} 
\left( u + \frac{(\alpha_i,\alpha_i)}{4}(m+1-2i) \right) \\
P_k(u) &=& 1, \mbox{ for } k\neq\l
\end{eqnarray*}
\end{defn}
We allow $m=0$, in which case $W_0(\l)$ is the trivial representation.

Viewed as a representation of $\g$, $\wml$ is a (not necessarily
irreducible) finite-dimensional representation in which the weight $m\wl$
occurs once and all other weights lie under $m\wl$ in the weight lattice.
The Kirillov-Reshetikhin formula deals specifically with a tensor product
of a number of such modules:
\begin{equation}
\label{def_decomp}
\bigotimes_{a=1}^N (W_{m_a}(\l_a) |_{\g})
\iso \Dsum_{\lambda} V_\lambda^{\dsum n_\lambda}
\end{equation}
where $V_\lambda$ is the irreducible $\g$-module with highest weight
$\lambda$ and it occurs $n_\lambda$ times in the decomposition of the
tensor product. 
Let us write $\wmax$ for $\sum_{a=1}^N m_a \w_{\l_a}$,
the highest weight (as a $\g$-module) of the tensor product.  Note that
$n_\wmax=1$. 

As discussed in chapter~\ref{chap_overview}, Kirillov and Reshetikhin
used the connections with mathematical physics to arrive at the following
conjecture of the multiplicities $n_\lambda$, where
$\lambda = \wmax - \sum n_i \alpha_i$.
\begin{equation}
\label{defZ}
n_\lambda = Z(\{\l\},\{m\}|n_1,\ldots,n_r) =
\sum_{\mbox{partitions}} \;\; \prod_{n\geq1} \;\; \prod_{k=1}^r
\binom{P^{(k)}_n(\nu) + \nu^{(k)}_n}{\nu^{(k)}_n}
\end{equation}
The sum is taken over all ways of choosing partitions
$\nu^{(1)},\ldots,\nu^{(r)}$ such that $\nu^{(i)}$ is a partition of
$n_i$ which has $\nu^{(i)}_n$ parts of size $n$ (so $n_i =
\sum_{n\geq1} n \nu^{(i)}_n$).  The function $P$ is defined by
\begin{eqnarray}
\label{defPgen}
P^{(k)}_n(\nu) &=& \sum_{a=1}^N \min(n,m_a)\delta_{k,\l_a}
  - 2 \sum_{h\geq 1} \min(n,h)\nu^{(k)}_{h} + \\
&&\hspace{1cm} +
  \sum_{j\neq k}^r \sum_{h\geq 1} \min(-c_{kj}n,-c_{j,k}h)\nu^{(j)}_{h}
\nonumber
\end{eqnarray}
where $C=(c_{ij})$ is the Cartan matrix of $\g$.  We define
${a\choose b}$ to be 0 whenever $a<b$.  Since the values of $P$ can be
negative, many of the binomial coefficients in~(\ref{defZ}) can be zero.

In Section~\ref{sec_algorithm}, we view the values of $P^{(k)}_n$ as the
coordinates of certain strings of weights of $\g$ which lie inside the
Weyl chamber.  This interpretation allows us to compute the values of
$n_\lambda$ much more efficiently.  Furthermore, the ``initial
substring'' relation on the labeling by strings of weights imposes the
structure of a rooted tree on the set of $\g$-modules which make up
$\wml$.

In section~\ref{sec_example}, we specialize to the case where $\g$ is
$\mathfrak{sl}_n$.  Here the representations of the Yangian are
irreducible when viewed as $\g$-modules, and the conjecture has already
been proven.  Therefore our results give a way to compute the
decomposition of a tensor product of representations of $\mathfrak{sl}_n$
whose associated Young diagrams are all rectangles.  The algorithm is
symmetric in all the factors, and again imposes the structure of a rooted
tree on the decomposition.

In Section~\ref{sec_growth}, we use this new tree structure to study the
asymptotics of the dimension of $\wml$ as $m$ gets large, based on the
fact that the tree structure of $\wml$ lifts to $W_{m+1}(\l)$.  We show
that the conjecture implies that the dimension grows asymptotically to a
polynomial in $m$, and compute the degree of this polynomial for every
simply-laced $\g$ and choice of $\wl$.

In Section~\ref{sec_table} we give a list of the decompositions of $\wml$
for all simply-laced $\g$ and small values of $m$ as derived numerically
from the conjecture, using the results of Section~\ref{sec_algorithm}.
For any choice of $\g$, representations $W_1(\l)$ are called fundamental
representations, since every finite-dimensional representation of $Y(\g)$
appears as a subquotient of a tensor product of such representations.
The decompositions of most of the fundamental representations were
calculated in~\cite{ChP} using completely different techniques, and those
calculations agree with ours.

\section{Structure in the simply-laced case}
\label{sec_algorithm}

Assume that our Lie algebra $\g$ of rank $r$ is simply-laced, and
otherwise retain the setup and notation of the previous section.  Then
equation~(\ref{defPgen}) simplifies to
\begin{equation}
\label{defP}
P^{(k)}_n(\nu) = \sum_{a=1}^N \min(n,m_a)\delta_{k,\l_a}
  - \sum_{j=1}^r c_{jk}
     \left( \sum_{h\geq 1} \min(n,h)\nu^{(j)}_{h} \right)
\end{equation}

Our goal is to find all choices for $\nu = (\nu^{(1)},\ldots,\nu^{(r)})$,
where $\nu^{(i)}$ is a partition of some integer $n_i$, such that
$P^{(k)}_n(\nu)$ is positive for all choices of $k$ and $n$.

\begin{thm}
\label{thm_decomp}
The pieces of the decomposition
$$
\bigotimes_{a=1}^N (W_{m_a}(\l_a) |_{\g})
\iso \Dsum_{\lambda} V_\lambda^{\dsum n_\lambda}
$$
arise from choices of partitions $\nu = (\nu^{(1)},\ldots,\nu^{(r)})$
which give a nonzero term in the sum in equation~(\ref{defZ}).  Such
choices are labeled by finite sequences $\d=(d_0,\ldots,d_s)$ of weights
of $\g$, with successive differences $\delta_i = d_i - d_{i-1}$ (and
$\delta_{s+1}=0$), such that:
  \begin{enumerate}
  \item[(i)] $d_0=0$ and $d_0\prec d_1\prec\cdots\prec d_s$,
  \item[(ii)] $\sum_{a=1}^N \min(n,m_a)\w_{\l_a} - d_n$ lies in the
    positive Weyl chamber for $0\leq n\leq s$, and
  \item[(iii)] $\delta_i \succeq \delta_{i+1}$ for all $1\leq i\leq s$.
  \end{enumerate}
where $\alpha\prec\beta$ means that $\beta-\alpha$ is in the cone of
positive roots of $\g$.  If we write $\wmax$ for $\sum_{a=1}^N m_a
\w_{\l_a}$, then the summand with label $\d=(d_0,\ldots,d_s)$ consists of
the $\g$-module of highest weight $\wmax - d_s$ with multiplicity
$$
\prod_{n\geq1} \;\; \prod_{k=1}^r \;
\binom{P^{(k)}_n(\d) + \d^{(k)}_n}{\d^{(k)}_n}
$$
where the values of $P^{(k)}_n(\d)$ and $\d^{(k)}_n$ are defined by the
relations
\begin{eqnarray*}
\sum_{a=1}^N \min(n,m_a)\w_{\l_a} - d_n &=& 
  \sum_{k=1}^r P^{(k)}_n(\d) \w_k \\
\delta_n - \delta_{n+1} &=& \sum_{k=1}^r \d^{(k)}_n \alpha_k
\end{eqnarray*}
All of these multiplicities are nonzero.
\end{thm}

\noindent {\bf Proof:}
Pick an arbitrary $\nu = (\nu^{(1)},\ldots,\nu^{(r)})$, where each
$\nu^{(i)}$ is a partition of some integer $n_i$.  Then for any
nonnegative integer $n$, the values $(P^{(1)}_n,\ldots,P^{(r)}_n)$ can be
thought of as the $\w$-coordinates of some weight; define
$$
\mu_n = \sum_{k=1}^r P^{(k)}_n \w_k
$$
A given $\nu$ contributes a nonzero term to the sum in~(\ref{defZ})
if and only if the corresponding weights $\mu_0=0, \mu_1, \mu_2,\ldots$ 
all lie in the dominant Weyl chamber.

The motivation for seeing these as weights is that the sum
in~(\ref{defP}) can be naturally realized as subtracting some linear
combination of roots.  If we let
\begin{equation}
\label{defd}
d_n = \sum_{k=1}^r
  \left( \sum_{h\geq 1} \min(n,h)\nu^{(k)}_{h} \right) \alpha_k
\end{equation}
then $\mu_n = \sum_{a=1}^N \min(n,m_a)\w_{\l_a} - d_n$.  Note that 
we have eliminated any reference to the Cartan matrix of $\g$.

Think of $\nu^{(1)},\ldots,\nu^{(r)}$ as Young diagrams, with 
$\nu^{(k)}$ having $\nu^{(k)}_{h}$ rows of length $h$.  The 
coefficient of $\alpha_k$ in $d_n$ is just the number of boxes in the 
first $n$ columns of $\nu^{(k)}$.  Now define $\delta_i = d_i - 
d_{i-1}$; if we write $\delta_n$ out as a linear combination of the 
roots $\{\alpha_i\}$, then the $\alpha_k$-coordinate is the number of 
boxes in the $n$th column of the Young diagram of $\nu^{(k)}$.  Thus a 
sequence of vectors $d_0=0, d_1, d_2,\ldots$ arises from partitions if 
and only if the $\delta_i$ are nonincreasing; that is, $\forall i\geq 
1 : \delta_i \succeq \delta_{i+1}$.

If we let $s$ be the size of the largest part in any of the partitions 
in $\nu$, then $d_s = d_t$ for all $t>s$ (and $s$ is the smallest 
index for which this is true).  So we label each summand of the 
decomposition with the (strictly increasing) sequence of weights $d_0=0 
\prec d_1\prec\cdots\prec d_s$.

We have constructed a label $\d=(d_0,\ldots,d_s)$ for each piece in 
the decomposition.  Conversely, from the label we can easily 
reconstruct the partitions, since the differences $\delta_i$ tell us 
the heights of successive columns of their Young diagrams.
\qed

We will say that the sequence $\d=(d_0,\ldots,d_s)$ has length $s$, and
we will call it {\em valid} if it satisfies conditions {\em (i)}, {\em
(ii)}, and {\em (iii)} above.  Note that the sequence of length 0
consisting of only $d_0=0$ is valid, arises from empty partitions, and
corresponds to the $V_{\wmax}$ component of the tensor product.

\begin{figure}
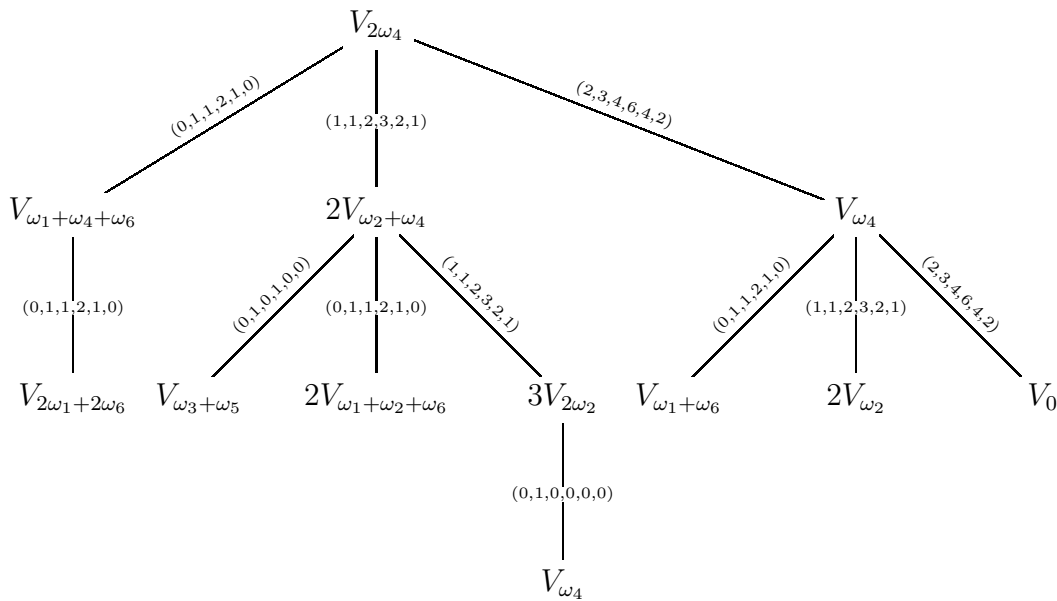

\begin{diagram}[labelstyle=\scriptscriptstyle]
&&&V_{2\w_4}&&&&&&& \\
&& \ldLine(3,2) ^{(0,1,1,2,1,0)} & \dLine ~{(1,1,2,3,2,1)} & 
    \rdLine(5,2) ^{(2,3,4,6,4,2)} &&&&&& \\
V_{\w_1 + \w_4 + \w_6}  & & &
    2 V_{\w_2 + \w_4} & & & & & V_{\w_4} & & \\
\dLine ~{(0,1,1,2,1,0)} & & \ldLine ^{(0,1,0,1,0,0)} &
    \dLine ~{(0,1,1,2,1,0)} & \rdLine ^{(1,1,2,3,2,1)} & & &
    \ldLine ^{(0,1,1,2,1,0)} & \dLine ~{(1,1,2,3,2,1)} &
    \rdLine ^{(2,3,4,6,4,2)} & \\
V_{2\w_1 + 2\w_6} & \;\; V_{\w_3 + \w_5} & &
    2 V_{\w_1 + \w_2 + \w_6} & & 3 V_{2\w_2} &
    \;\; V_{\w_1 + \w_6} & & 2 V_{\w_2} & & V_{0} \\
&&&&& \dLine ~{(0,1,0,0,0,0)} &&&&&\\
&&&&& V_{\w_4} &&&&&\\
\end{diagram}
\caption{Tree structure of the decomposition of $W_2(4)$ for $E_6$.}
\label{fig_tree}
\end{figure}

This decomposition is a refinement of the one in~(\ref{def_decomp}) 
since it is possible to find two different sequences $d_0,\ldots,d_s$ and 
$d'_0,\ldots,d'_t$ with $d_s = d'_t$.  This happens any time the sum 
in~(\ref{defP}) has more than one nonzero term.  One example of this 
occurs in $W_2(4)$ for $E_6$; see Figure~\ref{fig_tree}.

\begin{cor}
\label{cor_algor}
If $d_0,\ldots,d_s$ is a valid label then any initial segment 
$d_0,\ldots,d_{s'}$ (for $0 \leq s' < s$) is a valid label also.  
Conversely, given any valid label $d_0,\ldots,d_s$, we can extend it to 
another valid label by appending any weight $d_{s+1}$ which satisfies 
the conditions that $\min(s+1,m)\wl - d_{s+1}$ is in the 
positive Weyl chamber and, if $s>0$, that $d_s \prec d_{s+1} \preceq 
d_s+\delta_s$.
\end{cor}

\noindent
This follows immediately from conditions {\em (i)--(iii)}.
This is the key result which fails to hold true when $\g$ is not
simply-laced: pieces of the decomposition can still be given labels, but
it is possible to have an invalid label which can be extended to a valid
one.  Both the tree structure and the ease of generating labels are lost.
\qed

Since $d_0$ must be 0, this completely describes an effective 
algorithm for computing the conjectured decomposition.  The intuitive,
nondeterministic version is as follows:
\begin{algo}
\label{algo}
To decompose $\bigotimes_{a=1}^N (W_{m_a}(\l_a) |_{\g})$, let $\mu_0=0$ and
iterate the following steps for $n=1,2,3,\ldots$
\begin{enumerate}
\item
Add $\mbox{inc}_n$ to $\mu_{n-1}$, where $\mbox{inc}_n=\sum \w_{\l_a}$
for all $a$ such that $n\leq m_a$.
\item
Let $\mu_n$ be any weight in the Weyl chamber which is
$\mu_{n-1}+\mbox{inc}_n-\delta_n$, where $\delta_n$ is any sum of
positive roots less than or equal to $\delta_{n-1}$.  (If $n=1$, ignore
the $\delta_{n-1}$ part.)
\end{enumerate}
Stop when $\delta_n=0$.
\end{algo}
The computations in Section~\ref{sec_table} were computed using this
algorithm.  

Since truncating any label gives you another label, we can impose a 
tree structure on the parts of this decomposition, with a node of the 
tree corresponding to a summand in the decomposition from 
Theorem~\ref{thm_decomp}.  The ``children'' of the node with label 
$d_0,\ldots,d_s$ are all the nodes indicated by 
Corollary~\ref{cor_algor}; we can label the edges joining them to 
their parent with the various choices for the increment 
$\delta_{s+1}$.  For each $n\geq 0$, the $n$th row of the tree 
consists of all the nodes with labels of length $n$.

As an example of this structure, the tree for the decomposition of 
$W_2(4)$ for $\g=E_6$ is given in Figure~\ref{fig_tree}.  Scalars in 
front of modules, as in $2 V_{\w_2 + \w_4}$, indicate 
multiplicity.  The label $(a_1,\ldots,a_6)$ corresponds to an 
increment $\delta=\sum a_i\alpha_i$, so condition {\em (iii)} says 
that the labels along any path down from $V_{2\w_4}$ will be 
nonincreasing in each coordinate.  The labels on the edges are 
technically unnecessary, since they can be obtained by subtracting the 
highest weight of the child from the highest weight of the parent.  
However, as we will see in section~\ref{sec_tree}, they do record useful 
information that is not apparent by looking directly at the highest 
weights.

\section{Example: Tensor Products of Rectangles}
\label{sec_example} \label{sec_rectangles}

%
%
\newcommand{\mybox}{{\begin{picture}(15,15)%
\put(0,0){\line(1,0){15}}\put(0,0){\line(0,-1){15}}%
\put(15,0){\line(0,-1){15}}\put(0,-15){\line(1,0){15}}\end{picture}}}
\newcommand{\co}[2]{\multiput(#2,0)(0,-15){#1}{\mybox}}
\newcommand{\ron}[2]{\put(-3,-#2){\makebox(0,0)[r]{#1}}}
\newcommand{\bp}{\begin{picture}} \newcommand{\ep}{\end{picture}}
\newcommand{\rb}[1]{\raisebox{#1\unitlength}}

In this section we deal specifically with the situation when
$\g=\mathfrak{sl}_n$.  In this and only this case, every representation
of $\g$ is already a representation of $Y(\g)$: there is an evaluation
homomorphism from $Y(\g)$ to the universal enveloping algebra $U(\g)$,
which is the identity on the embedded copy $U(\g)\into Y(\g)$.  The
results of the previous section are therefore about decomposing tensor
products of representations of $\mathfrak{sl}_n$ whose associated Young
diagrams are rectangles.

In this case the decomposition formula was proved in~\cite{Ki} by
counting dimensions, and an explicit bijection between rigged
configurations and the Young tableaux that index pieces in the
decomposition was part of the detailed papers~\cite{KKR} and~\cite{KR1}
on the subject.  Thus in this section we give a purely combinatorial
algorithm for decomposing a tensor product of rectangles.

Note that the tensor product in equation~(\ref{def_decomp}) specializes
to a tensor of single columns when the coefficients $m_a$ are all 1.  In
this case the multiplicities $n_\lambda$ from equation~(\ref{defZ}) are
the Kostka numbers $K_{\lambda'\mu'}$, where $\lambda'$ is the partition
conjugate to $\lambda$ and $\mu'$ is the partition with columns of
heights $\l_1,\l_2,\ldots,\l_N$.

\subsection*{A Combinatorial Dictionary}

First we will fix our terminology for talking about Young diagrams, and
at the same time give the correspondence between the combinatorial view
and the one already presented.  For convenience, we will always assume we
are working with representations of $\mathfrak{sl}_n$ for sufficiently
large $n$; that is, for $n$ large enough that we never refer to
$\w_k$ for $k>n$.  (As a result, we could be working in
$\mathfrak{gl}_n$ just as easily as in $\mathfrak{sl}_n$.)

The representation $V_\lambda$ corresponds to a Young diagram, with each
fundamental weight $\w_k$ corresponding to a column of height $k$.
That is, if we write $\lambda = \sum_{k=1}^n b_k \w_k$, the Young
diagram has $b_k$ columns of height $k$.  The representations $\wml$ in
this case are just $V_{m\wl}$, corresponding to a rectangle with $m$
columns and $\l$ rows.

Given the Young diagram for $\lambda$, we get the diagram for
$\lambda-\alpha_i$ by moving one box down a row, from row $i$ to row
$i+1$.  Row $i+1$ must be shorter than row $i$, or else
$\lambda-\alpha_i$ is not inside the Weyl chamber.  The $\preceq$
relation of the previous section therefore translates into the {\em
dominance} ordering on Young diagrams: we say $Y_1$ dominates $Y_2$ if we
can obtain $Y_2$ from $Y_1$ by moving some boxes to lower rows of the
diagram.  Equivalently, $Y_1$ dominates $Y_2$ if the number of boxes in
the top $k$ rows of $Y_1$ is at least the number in the top $k$ rows of
$Y_2$, for every $k$.

The notion of moving down in the Weyl chamber is important enough that we
will invent a notation for recording how far down we have moved in the
dominance ordering.  We will record how many boxes move past the line
which is the bottom of the $k$th row; call these the {\em dominance
numbers} of the move.  We write the numbers in a column along the left
side of the diagram.  For example,
$$
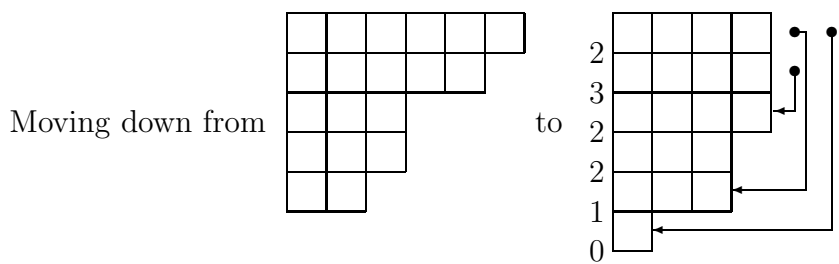
\mbox{Moving down from\,\, }
\rb{-45}{\bp(90,90)(0,-90)
\co5{0} \co5{15} \co4{30} \co2{45} \co2{60} \co1{75} \ep}
\mbox{ to }
\rb{-45}{\bp(100,90)(-15,-90)
\co6{0} \co5{15} \co5{30} \co3{45}
\ron2{15} \ron3{30} \ron2{45} \ron2{60} \ron1{75} \ron0{90}
\put(69,-7){\circle*{4}} \put(69,-22){\circle*{4}} \put(83,-7){\circle*{4}}
\put(69,-37){\line(0,1){15}} \put(69,-37){\vector(-1,0){8}}
\put(73,-67){\line(0,1){60}} \put(73,-67){\vector(-1,0){28}}
\put(73,-7){\line(-1,0){4}}
\put(83,-82){\line(0,1){75}} \put(83,-82){\vector(-1,0){68}} \ep}
\hspace{2cm}
$$
corresponds to subtracting
$2\alpha_1+3\alpha_2+2\alpha_3+2\alpha_4+\alpha_5$
from the weight
$\w_1+2\w_2+\w_4+2\w_5$.  The dominance numbers of the
move are $2,3,2,2,1$.

\begin{figure}
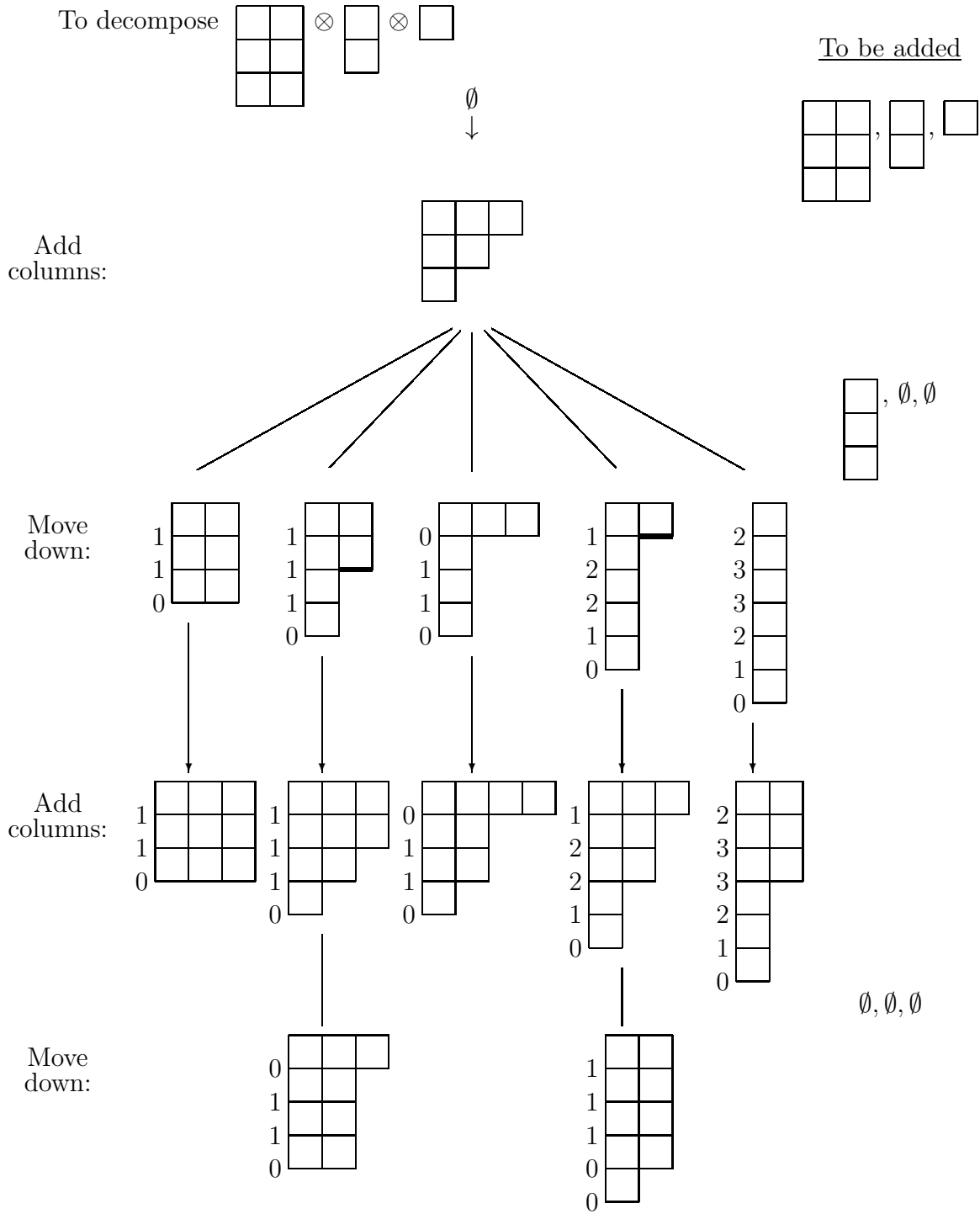

\begin{diagram}[w=5em]
\makebox(0,0)[l]{To decompose \  $
  \rb{-35}{\bp(30,45)(0,-45) \co3{0} \co3{15} \ep}\, \tensor \,
  \rb{-35}{\bp(15,45)(0,-45) \co2{0} \ep}\, \tensor \,
  \rb{-35}{\bp(15,45)(0,-45) \co1{0} \ep} $}
&&&&&& \mbox{\underline{To be added}} \\
&&& \mbox{\shortstack{$\emptyset$ \\ $\downarrow$}} &&& 
  \rb{-30}{\bp(30,45)(0,-45) \co3{0} \co3{15} \ep}\,,\,
  \rb{-30}{\bp(15,45)(0,-45) \co2{0} \ep}\,,\,
  \rb{-30}{\bp(15,45)(0,-45) \co1{0} \ep} \\
\raisebox{20pt}{\mbox{\shortstack{Add\\columns:}}} &&&
  {\bp(45,55)(0,-55) \co3{0} \co2{15} \co1{30} \ep} &&& \\
&& \ldLine(2,2)\ldLine(1,2) & \dLine & 
  \rdLine(1,2)\rdLine(2,2) &&
  \rb{-35}{\bp(15,45)(0,-45) \co3{0} \ep}\,,\,\emptyset,\emptyset \\
&\mbox{}&\mbox{}&\mbox{}&\mbox{}&\mbox{} \\
{\mbox{\shortstack{Move\\down:}}}
  & \rb{-25}{\bp(45,50)(-15,-50) 
     \co3{0} \co3{15} \ron1{15} \ron1{30} \ron0{45} \ep}
  & \rb{-40}{\bp(45,65)(-15,-65)
     \put(15,-30){\rule[-1\unitlength]{15\unitlength}{2\unitlength}}
     \co4{0} \co2{15} \ron1{15} \ron1{30} \ron1{45} \ron0{60} \ep}
  & \rb{-40}{\bp(60,65)(-15,-65) 
     \co4{0} \co1{15} \co1{30} \ron0{15} \ron1{30} \ron1{45} \ron0{60} \ep}
  & \rb{-55}{\bp(45,80)(-15,-80)
     \put(15,-15){\rule[-1\unitlength]{15\unitlength}{2\unitlength}}
     \co5{0} \co1{15} \ron1{15} \ron2{30} \ron2{45} \ron1{60} \ron0{75} \ep}
  & \rb{-70}{\bp(30,95)(-15,-95)
     \co6{0} \ron2{15} \ron3{30} \ron3{45} \ron2{60} \ron1{75} \ron0{90} \ep}
  \\
\makebox(0,20){}& \dTo & \dTo & \dTo & \dTo & \dTo \\
\mbox{\shortstack{Add\\columns:}}
  & \rb{-25}{\bp(60,50)(-15,-50) 
    \co3{0} \co3{15} \co3{30} \ron1{15} \ron1{30} \ron0{45} \ep}
  & \rb{-40}{\bp(60,65)(-15,-65)
    \co4{0} \co3{15} \co2{30} \ron1{15} \ron1{30} \ron1{45} \ron0{60} \ep}
  & \rb{-40}{\bp(75,65)(-15,-65)
    \co4{0} \co3{15} \co1{30} \co1{45}
      \ron0{15} \ron1{30} \ron1{45} \ron0{60} \ep}
  & \rb{-55}{\bp(60,80)(-15,-80) 
    \co5{0} \co3{15} \co1{30}
      \ron1{15} \ron2{30} \ron2{45} \ron1{60} \ron0{75} \ep}
  & \rb{-70}{\bp(45,95)(-15,-95)
    \co6{0} \co3{15}
      \ron2{15} \ron3{30} \ron3{45} \ron2{60} \ron1{75} \ron0{90} \ep}
  \\
&&\dLine&&\dLine&&\emptyset,\emptyset,\emptyset \\
\mbox{\shortstack{Move\\down:}}
  && \rb{-40}{\bp(60,65)(-15,-65)
     \co4{0} \co4{15} \co1{30} \ron0{15} \ron1{30} \ron1{45} \ron0{60} \ep}
  && \rb{-55}{\bp(45,80)(-15,-80) 
    \co5{0} \co4{15} \ron1{15} \ron1{30} \ron1{45} \ron0{60} \ron0{75} \ep}
  \\
\end{diagram}
\caption{Steps in decomposing a tensor of rectangles.}
\label{fig_algo}
\end{figure}

We can now rephrase Algorithm~\ref{algo} in terms of operations on Young
diagrams.  Figure~\ref{fig_algo} illustrates this process for a simple
example.

\begin{algo}
\label{algo_YD}
Every representation that occurs in the decomposition of the tensor
product of representations of $\mathfrak{sl}_n$ corresponding to
rectangular diagrams $R_1, R_2,\ldots,R_N$ can be found as follows.
Let $Y_0=\emptyset$, the empty Young diagram, and iterate the following
steps for $n=1,2,3,\ldots$
\begin{enumerate}
\item
Remove the first column from each of $R_1, R_2,\ldots,R_N$.  Add those
columns to $Y_{n-1}$ to form a new diagram, $Y'_n$.
\item
Let $Y_n$ be any diagram which can be obtained from $Y'_n$ by moving some
boxes of $Y'_n$ down to lower rows, such that the dominance numbers of
$Y_n$ are each less than or equal to the corresponding dominance number
for $Y_{n-1}$.  (Any numbers are allowed if $n=1$).
\end{enumerate}
Stop at any point.  (Equivalently, choose to not move down at all on the
next iteration, so all the dominance numbers are 0.)  There is a piece of
the decomposition corresponding to the final $Y_n$ with all remaining
columns of $R_1, R_2,\ldots,R_N$ added in.
\end{algo}

Finally, we need to address the multiplicities given by the binomial
coefficients in Theorem~\ref{thm_decomp}.  Suppose we carry out
Algorithm~\ref{algo_YD} and choose diagrams $Y_1,Y_2,\ldots,Y_s$ and then
stop.  The associated multiplicity is
$$
\prod_{n=1}^s \;\; \prod_{k=1}^r \;
\binom{P^{(k)}_n + \d^{(k)}_n}{\d^{(k)}_n}
$$
where this time the values of $P$ and $\d$ are given by
\begin{eqnarray*}
P^{(k)}_n &=& \mbox{number of columns of height $k$ in $Y_n$,} \\
\d^{(k)}_n &=&
  \mbox{decrease in $k$th completion number from $Y_n$ to $Y_{n+1}$}
\end{eqnarray*}
We consider all the completion numbers of $Y_{s+1}$ to be $0$.

\begin{figure}
\setlength{\unitlength}{.75\unitlength}
\begin{center}
\begin{diagram}[w=5em,h=1.3em]
&& \rb{-40}{\bp(60,55)(0,-55) \co3{0} \co3{15} \co2{30} \co1{45} \ep} \\
&& \mbox{} \\
& \ldLine(2,2)\ldLine(1,2) & \dLine & \rdLine(1,2)\rdLine(2,2) \\
\mbox{}&\mbox{}&\mbox{}&\mbox{}&\mbox{} \\
   \rb{-25}{\bp(45,50)(0,-50) \co3{0} \co3{15} \co3{30} \ep}
  & 2\,\rb{-40}{\bp(45,65)(0,-65) \co4{0} \co3{15} \co2{30} \ep}
  & \rb{-40}{\bp(60,65)(0,-65) \co4{0} \co3{15} \co1{30} \co1{45} \ep}
  & 2\,\rb{-55}{\bp(45,80)(0,-80) \co5{0} \co3{15} \co1{30} \ep}
  & \rb{-70}{\bp(30,95)(0,-95) \co6{0} \co3{15} \ep}
  \\
& \dLine && \dLine \\
\mbox{} \\
&  \rb{-40}{\bp(45,65)(0,-65) \co4{0} \co4{15} \co1{30} \ep}
 && \rb{-55}{\bp(30,80)(0,-80) \co5{0} \co4{15} \ep}
\end{diagram}
\end{center}
\vspace{-1em}
\caption{Decomposition of
$V_{2\w_3}\tensor V_{\w_2}\tensor V_{\w_1}$
(with multiplicities).}
\label{fig_rect}
\end{figure}

Graphically, this means we get a contribution to the multiplicity when we
move from $Y_n$ to $Y_{n+1}$ if the $k$th completion number decreases and
the $k$th row of $Y_n$ overhangs the $k+1$st row.  If $n\neq s$ then the
multiplicities of the entire section of the tree from $Y_{k+1}$ down have
that binomial coefficient as a factor.  Unfortunately, this never happens
in Figure~\ref{fig_algo}; the smallest examples in which this comes up
are too large to include here.  (The interested reader can see it take
place by trying the example in Figure~\ref{fig_algo} with $2\w_2$
instead of $\w_2$.)

If $n=s$ then we are looking for an instance where the $k$th row of $Y_n$
overhangs the $k+1$st next to any completion number which is nonzero; in
this case the factor affects the multiplicity of only the $Y_s$ node on
the tree, not those below it.  This happens twice in Figure~\ref{fig_algo};
the overhangings are marked in bold.  The complete decomposition,
including multiplicities, is presented in Figure~\ref{fig_rect}.

\section{The Growth of Trees}
\label{sec_tree} \label{sec_growth}

In this section we stop talking about tensor products and just look at 
decompositions of the modules $\wml$ themselves.  We show that the 
trees of the decompositions of $\wml$ for different values of $m$ are 
compatible with one another.  As an application of this newfound 
structure, we will prove that as $m$ gets large, the dimension of the 
representation $\wml$ grows like a polynomial in $m$, and will give a 
method to compute the degree of the polynomial growth.  All statements 
assume the conjectural formulas for multiplicities of $\g$-modules.
Roots and weights are numbered as in Table~\ref{table_dd}
(p.~\pageref{table_dd}).

We begin with another corollary of Theorem~\ref{thm_decomp}, whose 
notation we retain.

\begin{cor}
\label{cor_lift}
If $d_0,\ldots,d_s$ is a valid label for a piece of the decomposition 
of $W_m(\l)$, then it is also a valid label for $W_{m'}(\l)$ for any 
$m'>m$, and for any $m'\geq s$.
\end{cor}

\noindent Both parts are based on the fact that condition {\em (ii)} 
of Theorem~\ref{thm_decomp} is the only one that depends on $m$.  For 
$m'>m$, if $\min(m,n)\wl - d_n$ is a nonnegative linear combination of 
the $\{\wl\}$ then adding some nonnegative multiple of $\wl$ will not 
change that fact.  For $m'\geq s$, the value of $m'$ is 
irrelevant; the weights we look at are just $n\wl - d_n$ for $0\leq 
n\leq s$.  \qed

If we can lift labels from $W_m(\l)$ to $W_{m+1}(\l)$, we can also lift
the entire tree structure.  Specifically, the lifting of labels extends
to a map from the tree of $W_m(\l)$ to the tree of $W_{m+1}(\l)$ which
preserves the increment $\delta$ of each edge and lifts each $V_\lambda$
to $V_{\lambda+\wl}$.  The $m'\geq s$ part of 
Corollary~\ref{cor_lift} tells us that this map is a bijection on rows
$0,1,\ldots,m$ of the trees, where the labels have length $s\leq m$.  On
this part of the tree, multiplicities are also preserved.  This 
follows from the formula for multiplicities in Theorem~\ref{thm_decomp}:
the only values of $P^{(k)}_n(\d)$ that change are for $n=m+1$, but
$\d^{(k)}_n=0$ when $n$ is greater than the length of the label, so the
product of binomial coefficients is unchanged.

\begin{defn}
Let $\tl$ be the tree whose top $n$ rows are those of $\wml$ for all 
$m\geq n$.
\end{defn}

The highest weight associated with an individual node appearing in 
$\tl$ is only well-defined up to addition of any multiple of $\wl$, 
but the difference $\delta$ between any node and its parent is 
well-defined.  (These differences are the labels on the edges of the 
tree in Figure~\ref{fig_tree}.)  We can characterize each node by the 
string of successive differences $\delta_1 \succeq \delta_2 
\succeq\cdots\succeq \delta_s$ which label the $s$ edges in the path 
from the root of the tree to that node.  The multiplicity of a node of 
$\tl$ is well-defined, as already noted.

\bigskip \noindent
{\bf Claim: } {\em For a fixed $\l$, the dimension of $\wml$ grows as a
polynomial in $m$, whose degree we can calculate.}
\bigskip

We will study the growth of the tree $\tl$.  The precise statement of the
claim is in Theorem~\ref{thm_growth}.

The tree of $\wml$ matches $\tl$ exactly in the top $m$ rows.  The number
of rows in the tree of $\wml$ is bounded by the largest
$\alpha$-coordinate of $m\wl$, since if $\delta_1,\ldots,\delta_s$
is a label of $\wml$ then $m\wl - \sum_{i=1}^s \delta_i$ must be
in the positive Weyl chamber, and whatever $\alpha$-coordinate is nonzero
in $\delta_s$ must be nonzero in all of the $\delta_i$.  Therefore to
prove that the dimension of $\wml$ grows as a polynomial in $m$, it
suffices to prove that the dimension of the part of $\wml$ which
corresponds to the top $m$ rows of $\tl$ does so.

Now we need to examine the structure of the tree $\tl$.  The path 
$\delta_1,\ldots,\delta_s$ to reach a vertex is a sequence of weights 
whose $\alpha$-coordinates are nonincreasing.  Write this instead as 
$\Dtmt$ where the $\Delta_i$ are strictly decreasing and $m_i$ is the 
number of times $\Delta_i$ occurs among $\delta_1,\ldots,\delta_s$; we 
will say this path has {\em path-type} $\Dt$.  The number of 
path-types that can possibly appear in the tree $\tl$ is finite, since 
each $\Delta_i$ is between $\wl$ and 0 and has integer 
$\alpha$-coordinates.

We need to understand which path-types $\Dt$ and which choices of
exponents $m_i$ correspond to paths which actually appear in $\tl$.
Given a path $\delta_1,\ldots,\delta_s$, assume that $m>s$ and recall
$\mu_n = n\wl - d_n = n\wl - \sum_{i=1}^n \delta_i$.
Condition~{\em (ii)} from Theorem~\ref{thm_decomp} requires that $\mu_n$
is in the positive Weyl chamber for $1\leq n\leq s$; that is, the
$\w$-coordinates of $\mu_n$ must always be nonnegative.  (These
coordinates are just the values of $P^{(k)}_n$ from
Theorem~\ref{thm_decomp}.)  Since $\mu_n = \mu_{n-1} + \wl -
\delta_n$, we need to keep track of which $\w$-coordinates of
$\wl-\delta_n$ are positive and which are negative.

\begin{defn}
For a path-type $\Dt$, we say that $\Delta_i$ {\em provides} 
$\w_k$ if the $\w_k$-coordinate of $\wl - \Delta_i$ is 
positive, and that it {\em requires} $\w_k$ if the coordinate is 
negative.
\end{defn}

Geometrically, $\Delta_i$ providing $\w_k$ means that each $\Delta_i$
in the path moves the sequence of $\mu$s away from the $\w_k$-wall of
the Weyl chamber, while requiring $\w_k$ moves towards that wall.
The terminology is justified by restating what condition {\em (ii)}
implies about path-types in these terms:

\begin{lemma}
\label{lem_pathtypes}
The tree $\tl$ contains paths of type $\Dt$ if and only if, for every 
$\Delta_n$, $1\leq n\leq t$, every $\w_i$ required by $\Delta_n$ 
is provided by some $\Delta_k$ with $k<n$.
\end{lemma}

\noindent
The ``only if'' part of the equivalence is immediate from the preceding
discussion: the sequence $\mu_0,\mu_1,\ldots$ starts at $\mu_0=0$, and if
it moves towards any wall of the Weyl chamber before first moving 
away from it, it will pass through the wall and some $\mu_i$ will be
outside the chamber.  Conversely, if $\Dt$ is any path-type which
satisfies the condition of the lemma, then $\Dtmt$ will definitely appear
in the tree when $m_1 \gg m_2 \gg\cdots\gg m_t$.  This ensures that the
coordinates of the $\mu_i$ are always nonnegative, since the sequence of
$\mu$s moves sufficiently far away from any wall of the Weyl chamber
before the first time it moves back towards it.
\qed

We could compute the exact conditions on the $m_i$ for a specific path;
in general, they all require that $m_n$ be bounded by some linear
combination of $m_1,\ldots,m_{n-1}$, and the first $m_i$ appearing with
nonzero coefficient in that linear combination has positive coefficient.

Now we can show that the number of nodes of path-type $\Dt$ appearing on
the $m$th level of the tree grows as $m^{t-1}$.  Consider the path
$\Dtmt$ as a point $(m_1,\ldots,m_t)$ in $\R^t$.  The path ends on row
$m$ if $m=m_1+\cdots+m_t$, so solutions lie on a plane of dimension
$t-1$; the number of solutions to that equality in nonnegative integers
is ${m+t-1}\choose{t-1}$, which certainly grows as $m^{t-1}$, as
expected.  The further linear inequalities on the $m_i$ which ensure that
$\mu_1,\ldots,\mu_m$ remain in the Weyl chamber correspond to hyperplanes
through the origin which our solutions must lie on one side of, but the
resulting region still has full dimension $t-1$ since the generic point
with $m_1 \gg m_2 \gg\cdots\gg m_t$ satisfies all of the inequalities, as
shown above.

The highest weight of the $\g$-module at the node associated with the
generic solution of the form $m_1 \gg m_2 \gg\cdots\gg m_t$ grows
linearly in $m$.  Its dimension, therefore, grows as a polynomial in $m$,
and the degree of the polynomial is just the number of positive roots of
the Lie algebra which are not orthogonal to the highest weight.  The only
positive roots perpendicular to this generic highest weight are those
perpendicular to every highest weight which comes from a path of type
$\Dt$, and the number of such roots is the degree of polynomial growth of
the dimensions of the representations of the $\g$-module.

We can figure out how the multiplicities of nodes with a specific 
path-type grow as well.  Theorem~\ref{thm_decomp} gives a formula for 
multiplicities as a product of binomial coefficients over $1\leq k\leq 
r$ and $n\geq 1$.  The only terms in the product which are not 1 
correspond to nonzero values of $\delta_n - \delta_{n+1}$.  In the 
path $\Dtmt$, these occur only when $n=m_1+\cdots+m_i$ for some $1\leq 
i\leq t$, so that $\delta_n - \delta_{n+1}$ is $\Delta_i - 
\Delta_{i+1}$ (where $\Delta_{t+1}$ is just 0).  Following our 
previous notation, let $\delta_n - \delta_{n+1} = \d_n = \sum 
\d^{(k)}_n \alpha_k$.  If we take any $k$ for which $\d^{(k)}_n$ is 
nonzero, there are two possibilities for the contribution to the 
multiplicity from its binomial coefficient.  If $\w_k$ has been 
provided by at least one of $\Delta_1,\ldots,\Delta_i$, then the value 
$P^{(k)}_n$ is a linear combination of $m_1,\ldots,m_i$, which grows 
linearly as $m$ gets large.  In this case, the binomial coefficient 
grows as a polynomial in $m$ of degree $\d^{(k)}_n$.  On the other 
hand, if $\w_k$ has not been provided, then the binomial 
coefficient is just 1.

For any $k$, $1\leq k\leq r$, define $f(k)$ to be the smallest $i$ in 
our path-type such that $\Delta_i$ provides $\w_k$; we say that 
$\Delta_i$ provides $\w_k$ for the first time.  Then the total 
contribution to the multiplicity from the coordinate $k$ will be the 
product of the contributions when $n=m_1+\cdots+m_j$ for $j=f(k), 
f(k)+1, \ldots, t$.  As $m$ gets large, the product of these 
contributions grows as a polynomial of degree $\sum_{j=f(k)}^t 
\d^{(m_1+\cdots+m_j)}_n$; that is, the sum of the decreases in the 
$\alpha_k$-coordinate of the $\Delta$s.  But since $\Delta_{t+1}$ is 
just 0, that sum is exactly the $\alpha_k$-coordinate of 
$\Delta_{f(k)}$.

So given a path-type $\Dt$ which Lemma~\ref{lem_pathtypes} says appears
in $\tl$, the total of the multiplicities of the nodes of that path-type
which appear in the top $m$ rows of $\tl$ grows as a polynomial of degree
\begin{equation}
\label{defg}
g(\Dt) = t + \sum_{k=1}^{r} \alpha_k\mbox{-coordinate of }\Delta_{f(k)}
\end{equation}
where we take $\Delta_{f(k)}$ to be 0 if $\w_k$ is not provided by any
$\Delta$ in the path-type.  This value is just the sum of the degrees of
the polynomial growths described above.

Finally, since there are only finitely many path-types, the growth of the
entire tree $\tl$ is the same as the growth of the part corresponding to
any path-type $\Dt$ which maximizes $g(\Dt)$.  So we have proven the
following, up to some calculation:

\begin{thm}
\label{thm_growth}
Let $\g$ be simply-laced with decompositions of $\wml$ given by
Theorem~\ref{thm_decomp}.  Then the dimension of the representation
$\wml$ as $m$ gets large is asymptotic to a polynomial in $m$ of degree
$\mathrm{perp}(\Dt)+g(\Dt)$, where the path-type $\Dt$ is one which
maximizes the value of $g$, and $\mathrm{perp}(\Dt)$ is the number of
positive roots orthogonal to all highest weights of nodes with path-type
$\Dt$.

\begin{enumerate}
\item
If $\g$ is of type $A_n$ then the maximum value of $g(\Dt)$ is $0$, for 
all $1\leq\l\leq n$.
\item
If $\g$ is of type $D_n$ then the maximum value of $g(\Dt)$ is 
$\lfloor\l/2\rfloor$, for $1\leq\l\leq n-2$, and $0$ for $\l=n-1,n$.
\item
If $\g$ is of type $E_6$, $E_7$, or $E_8$, the maximum value of 
$g(\Dt)$ is
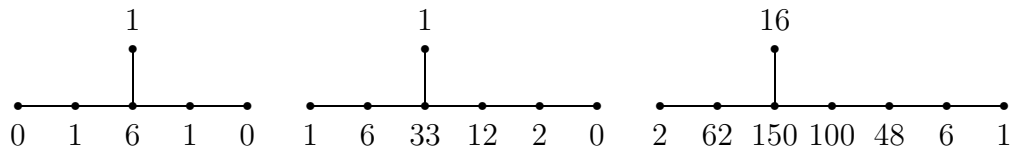
\begin{center} \setlength{\unitlength}{.3in}
\newcommand{\num}[1]{{\makebox(0,1){${#1}$}}} 
\newcommand{\vx}{{\circle*{.15}}}
\begin{picture}(4,3) 
\multiput(0,1)(1,0){5}{\vx}\put(2,2){\vx}
\put(0,1){\line(1,0){4}}\put(2,1){\line(0,1){1}}
\put(2,2){\num{1}}\put(0,0){\num{0}}\put(1,0){\num{1}}
\put(2,0){\num{6}}\put(3,0){\num{1}}\put(4,0){\num{0}}
\end{picture}\hfill
\begin{picture}(5,3) 
\put(2,2){\vx}\multiput(0,1)(1,0){6}{\vx}
\put(0,1){\line(1,0){5}}\put(2,1){\line(0,1){1}}
\put(2,2){\num{1}}\put(0,0){\num{1}}\put(1,0){\num{6}}\put(2,0){\num{33}}
\put(3,0){\num{12}}\put(4,0){\num{2}}\put(5,0){\num{0}}
\end{picture}\hfill
\begin{picture}(6,3) 
\put(2,2){\vx}\multiput(0,1)(1,0){7}{\vx}
\put(0,1){\line(1,0){6}}\put(2,1){\line(0,1){1}}
\put(2,2){\num{16}}\put(0,0){\num{2}}\put(1,0){\num{62}}\put(2,0){\num{150}}
\put(3,0){\num{100}}\put(4,0){\num{48}}\put(5,0){\num{6}}\put(6,0){\num{1}}
\end{picture}\hspace*{\fill}
\end{center}
\end{enumerate}
\end{thm}

We will complete the proof by exhibiting the path-types which give 
the indicated values of $g$ and proving they are maximal.

If $\Dt$ maximizes the value of $g$, then it cannot be obtained from any
other path-type by inserting an extra $\Delta$, since any insertion
would increase the length $t$ and would not decrease the sum in the
definition of $g$.  Therefore each $\Delta_k$ in our desired path-type
must be in the positive root lattice, allowable according to
Lemma~\ref{lem_pathtypes}, and must be maximal (under $\preceq$) in
meeting those requirements; we will call a path-type maximal if this is
the case.

In particular, if $\wl$ is in the root lattice then $\Delta_1$ will
be $\wl$, and a $g$-value of 0 corresponds exactly to an $\wl$ which is
not in the root lattice and is a minimal weight.  Thus the 0s above
can be verified by inspection; these are exactly the cases in which
$\wml$ remains irreducible as a $\g$-module.  Similarly, if $\wl$ is not
in the root lattice but there is only one point in the lattice and in
the Weyl chamber under $\wl$, the path-type will consist just of that
point.  We can now limit ourselves to path-types of length greater than
one.

If $\g$ is of type $D_n$ then for each $\wl$, $2\leq\l\leq n-2$, 
there is a unique maximal path-type:
$$
\begin{array}{ll}
\wl \succ \wl-\w_2 \succ \wl-\w_4 \succ \cdots \succ \wl-\w_{\l-2}
  & \mbox{when $\l$ is even} \\
\wl-\w_1 \succ \wl-\w_3 \succ \cdots \succ \wl-\w_{\l-2}
  & \mbox{when $\l$ is odd}
\end{array}
$$
In both cases, the only contribution to $g$ comes from the length of 
the path, which is $\lfloor\l/2\rfloor$.  This also means that the 
nodes of the tree $\tl$ will all have multiplicity 1 in this case.  
See figure~\ref{fig_Dtree} (p.~\pageref{fig_Dtree}) for a graphical 
interpretation of this path-type.

When $\g$ is of type $E_6$, $E_7$ or $E_8$, the following weights have a
unique maximal path-type (of length $>1$), whose $g$-value is given in
Theorem~\ref{thm_growth}:
$$
\begin{array}{lll}
E_6
& \l=4 & \w_4 \succ \w_4-\w_2 \succ \w_4-\w_1-\w_6
       \succ \w_2+\w_4-\w_3-\w_5 \succ 2\w_2-\w_4 \\
E_7
& \l=3 & \w_3 \succ \w_3-\w_1 \succ \w_3-\w_6 \succ
       \w_1+\w_6-\w_4 \succ 2\w_1-\w_3 \\
& \l=6 & \w_6 \succ \w_6-\w_1 \\
E_8
& \l=1 & \w_1 \succ \w_1-\w_8 \\
& \l=7 & \w_7 \succ \w_7-\w_8 \succ \w_7-\w_1 \succ
       \w_7+\w_8-\w_6 \succ 2\w_8-\w_7 \\
& \l=8 & \w_8
\end{array}
$$
We consider the remaining weights in $E_8$ next.  Consider the
incomplete path-type
$$
\wl \succ \wl-\w_8 \succ \wl-\w_1 \succ \wl-\w_6+\w_8
\succ \wl+\w_1-\w_4+\w_8 \succ \cdots
$$
where $\wl$ is any fundamental weight which is in the root lattice 
and high enough that all of the weights in question lie in the Weyl 
chamber.  The path so far provides $\w_8$, $\w_1$, $\w_6$ 
and $\w_4$; notice that for any $\w_i$ which has not been 
provided, all of its neighbors in the Dynkin diagram have.  Therefore 
we can extend this path four more steps by subtracting one of 
$\alpha_2$, $\alpha_3$, $\alpha_5$ and $\alpha_7$ at each step, to 
produce a path in which every $\w_i$ has been provided.  This can 
be extended to a full path-type by subtracting any $\alpha_i$ at each 
stage until we reach the walls of the Weyl chamber.

The resulting path-type is maximal, and is the unique maximal one up to a
sequence of transformations of the form
$$
\cdots \succ \Delta \succ \Delta-\lambda \succ \Delta-\lambda-\mu \succ \cdots
\mapsto
\cdots \succ \Delta \succ \Delta-\mu \succ \Delta-\lambda-\mu \succ \cdots
$$
which do not affect the rate of growth $g$.  All relevant weights are in
the Weyl chamber if and only if $\wl\succ\xi=(4,8,10,14,12,8,6,2)$; this
turns out to be everything except $\w_1$, $\w_7$ and $\w_8$,
whose path-types are given above.  If the path-type could start at $\xi$,
it would have growth $g=8$, though this is not possible since the last
weight in the path-type would be $0$ in this case.  But each increase of
the starting point of the path by any $\alpha_i$ increases $g$ by 2
(1 from the length of the path and 1 from the multiplicity).  So the growth
for any $\wl\succ\xi$ is a linear function of its height with coefficient
2; $g=2\hyt(\wl)-120$.

The only remaining cases are $\w_4$ and $\w_5$ when $\g$ is of
type $E_7$.  Both work like the general case for $E_8$, beginning instead
with the incomplete path-types
$$
\begin{array}{ll}
\l=4 & \w_4 \succ \w_4-\w_1 \succ \w_4-\w_6
     \succ \w_1 \succ \cdots \\
\l=5 & \w_5-\w_7 \succ \w_5-\w_2 \succ
     \w_5+\w_7-\w_1-\w_6 \succ
     \w_2+\w_7-\w_3 \succ \cdots
\end{array}
$$
This concludes the proof of Theorem~\ref{thm_growth}.
\qed

The same argument used for $E_8$ shows that for any choice of $\g$, all
``sufficiently large'' weights $\wl$ in a particular translate of the
root lattice will have growth given by $2\hyt(\wl)-c$ for some fixed $c$.
A weight is sufficiently large if every $\w_i$ is provided in its maximal
path.  Thus we can easily check that $\w_4$ and $\w_5$ qualify for $E_7$,
and in both cases $c=63$.  Similarly, $\w_4$ for $E_6$ qualifies, and
$c=36$.  While there are no sufficiently large fundamental weights for
$A_n$ or $D_n$, we can compute what the maximal path-type would be if one
did exist, and in all cases, $c$ is the number of positive roots.  A
uniform explanation would be nice, though the exhaustive computation does
provide a complete proof.

\section{Computations}
\label{sec_table}

This section gives the decompositions of $\wml$ into $\g$-modules
predicted by the conjectural formulas in~\cite{KR} for $D_n$ and $E_n$.
We also give the tree structure defined in Section~\ref{sec_algorithm}.
Roots and weights are numbered as in Table~\ref{table_dd}
(p.~\pageref{table_dd}).  As already noted, when $\g$ is of type $A_n$,
the $Y(\g)$-modules $\wml$ remain irreducible when viewed as
$\g$-modules.

The representations $\wml$ when $m=1$ are called fundamental
representations.  In the setting of $U_q(\g)$-module decompositions of
$U_q(\ghat)$ modules, the decompositions of the fundamental
representations for all $\g$ and most choices of $\wl$ appear
in~\cite{ChP}, calculated using techniques unrelated to the conjecture
used in~\cite{KR} to give formulas~(\ref{defZ}) and~(\ref{defPgen}).
Those computations agree with the ones given below.

The choices of $\wl$ not calculated in~\cite{ChP} are exactly those in
which the maximal path-type (Theorem~\ref{thm_growth}) is not unique.

\subsection*{\boldmath $D_n$}

Let $\g$ be of type $D_n$.  As already noted, the fundamental weights
$\w_{n-1}$ and $\w_{n}$ are minimal with respect to $\preceq$,
so $W_m(n-1)$ and $W_m(n)$ remain irreducible as $\g$-modules.
Now suppose $\l\leq n-2$.  Then the structure of the weights in the
Weyl chamber under $\w_\l$ does not depend on $n$, and so the
decomposition of $\wml$ in $D_n$ is the same for any $n \geq \l+2$.

\begin{figure}
\newcommand{\pt}{\bullet}
%
%
\begin{diagram}[width=2em,height=1.7em,abut]
&&&\raisebox{-6pt}{\Large $V_{3\w_6}$} \\
&&&\pt \\
&&\ldLine(1,2)&&\rdLine(1,2)\rdLine(6,2) \\
&&\pt&&\pt&&&  &&\pt&& \\ 
&\ldLine(1,2)&&\ldLine(1,2)&&\rdLine(1,2)&& 
  &\ldLine(1,2)&&\rdLine(1,2)\rdLine(4,2) \\
&\pt&&\pt&&\pt&&  &\pt&&\pt& &&\pt \\
\ldLine(1,2)&&\ldLine(1,2)&&\ldLine(1,2)&&\rdLine(1,2)&
  \ldLine(1,2)&&\ldLine(1,2)&&\rdLine(1,2)&
  \ldLine(1,2)&&\rdLine(1,2)\rdLine(3,2) \\
\pt&&\pt&&\pt&&\pt & \pt&&\pt&&\pt & \pt&&\pt && \pt
\end{diagram}

\hfill
\parbox[t]{2.75in}{\renewcommand{\baselinestretch}{.95}\small
The children of a node with highest weight $\lambda$ have highest
weights $\lambda-\w_6+\w_4$, $\lambda-\w_6+\w_2$, and
$\lambda-\w_6$, if they are to the left, right, or far right of the
parent, respectively.}
\hfill
\parbox[t]{2.75in}{\renewcommand{\baselinestretch}{.95}\small
The tree for $W_3(7)$ looks identical, but with $V_{3\w_7}$ on top
and highest weights $\lambda-\w_7+\w_5$,
$\lambda-\w_7+\w_3$, and $\lambda-\w_5+\w_1$ for the
children.}
\hfill
\caption{Tree structure of the decomposition of $W_3(6)$ for $D_n$ for
any $n\geq 8$.}
\label{fig_Dtree}
\end{figure}

As mentioned in the proof of Theorem~\ref{thm_growth}, there is a
unique maximal path-type for each $\wl$, and there are no multiplicities
greater than 1.  The decomposition is therefore very simple: if
$\l\leq n-2$ is even, then
$$
\wml \iso \Dsum_{k_2+k_4+\ldots+k_{\l-2}+k_\l = k \leq m}
     V_{ k_2 \w_2 + k_4 \w_4 + \ldots + k_{\l-2} \w_{\l-2}
         + (m-k)\w_\l }
$$
and if $\l$ is odd, then
$$
\wml \iso \Dsum_{k_1+k_3+\ldots+k_{\l-2} = k \leq m}
     V_{ k_1 \w_1 + k_3 \w_3 + \ldots + k_{\l-2} \w_{\l-2}
         + (m-k)\w_\l }
$$
The minor difference is because $\w_\l$ for $\l$ odd is not in
the root lattice.  The sum $k$ is the level of the tree on which that
module appears, and the parent of a module is obtained by subtracting 
1 from the first of $k_{\l-2}, k_{\l-4},\ldots$ which is nonzero (or from
$k_\l$ if nothing else is nonzero and $\l$ is even).

Figure~\ref{fig_Dtree} illustrates the tree structure for $\wml$ 
when $\l=6$ or $\l=7$ (and $n\geq\l+2$).  For $\l=4$ or $5$, the tree 
would look like the left-most triangle of figure~\ref{fig_Dtree}, 
while for $8$ or $9$ the evident recursive pattern would be one level 
deeper.

\subsection*{\boldmath $E_n$}

When $\g$ is of type $E_n$ the tree structure is much more irregular:
these are the only cases in which a $\g$-module can appear in more that
one place in the tree and in which a node on the tree can have
multiplicity greater than one.

We indicate the tree structure as follows: we list every node in the
tree, starting with the root and in depth-first order, and a node on
level $k$ of the tree is written as $\oplustag{k} V_\lambda$.  This is
enough information to recover the entire tree, since the parent of
that node is the most recent summand of the form $\oplustag{k-1} V_\mu$.
Comparing Figure~\ref{fig_tree} to its representation here should make
the notation clear.

Due to space considerations, for $E_6$ we list calculations for $m\leq 
3$, for $E_7$ we list $m\leq 2$, and for $E_8$ only $m=1$.  The 
tree decomposition for $W_3(4)$ for $E_7$, for example, would have 
836 components.

\begin{enumerate} \sloppy
\item[{\boldmath $E_6$}]
\item[$W_m(1)$] remains irreducible for all $m$.
\item[$W_1(2)$] $\iso V_{\w_2}
 \oplustag{1} V_{0}$
\item[$W_2(2)$] $\iso V_{2 \w_2}
 \oplustag{1} V_{\w_2}
 \oplustag{2} V_{0}$
\item[$W_3(2)$] $\iso V_{3 \w_2}
 \oplustag{1} V_{2 \w_2}
 \oplustag{2} V_{\w_2}
 \oplustag{3} V_{0}$
\item[$W_1(3)$] $\iso V_{\w_3}
 \oplustag{1} V_{\w_6}$
\item[$W_2(3)$] $\iso V_{2 \w_3}
 \oplustag{1} V_{\w_3+\w_6}
 \oplustag{2} V_{2 \w_6}$
\item[$W_3(3)$] $\iso V_{3 \w_3}
 \oplustag{1} V_{2 \w_3+\w_6}
 \oplustag{2} V_{\w_3+2 \w_6}
 \oplustag{3} V_{3 \w_6}$
\item[$W_1(4)$] $\iso V_{\w_4}
 \oplustag{1} V_{\w_1+\w_6}
 \oplustag{1} 2 V_{\w_2}
 \oplustag{1} V_{0}$
\item[$W_2(4)$] $\iso V_{2 \w_4}
 \oplustag{1} V_{\w_1+\w_4+\w_6}
 \oplustag{2} V_{2 \w_1+2 \w_6}
 \oplustag{1} 2 V_{\w_2+\w_4}
 \oplustag{2} V_{\w_3+\w_5}
 \oplustag{2} 2 V_{\w_1+\w_2+\w_6}
 \oplustag{2} 3 V_{2 \w_2}
 \oplustag{3} V_{\w_4}
 \oplustag{1} V_{\w_4}
 \oplustag{2} V_{\w_1+\w_6}
 \oplustag{2} 2 V_{\w_2}
 \oplustag{2} V_{0}$
\item[$W_3(4)$] $\iso V_{3 \w_4}
 \oplustag{1} V_{\w_1+2 \w_4+\w_6}
 \oplustag{2} V_{2 \w_1+\w_4+2 \w_6}
 \oplustag{3} V_{3 \w_1+3 \w_6}
 \oplustag{1} 2 V_{\w_2+2 \w_4}
 \oplustag{2} V_{\w_3+\w_4+\w_5}
 \oplustag{2} 2 V_{\w_1+\w_2+\w_4+\w_6}
 \oplustag{3} V_{\w_1+\w_3+\w_5+\w_6}
 \oplustag{3} 2 V_{2 \w_1+\w_2+2 \w_6}
 \oplustag{2} 3 V_{2 \w_2+\w_4}
 \oplustag{3} V_{2 \w_4}
 \oplustag{3} 2 V_{\w_2+\w_3+\w_5}
 \oplustag{3} 3 V_{\w_1+2 \w_2+\w_6}
 \oplustag{4} V_{\w_1+\w_4+\w_6}
 \oplustag{3} 4 V_{3 \w_2}
 \oplustag{4} 2 V_{\w_2+\w_4}
 \oplustag{1} V_{2 \w_4}
 \oplustag{2} V_{\w_1+\w_4+\w_6}
 \oplustag{3} V_{2 \w_1+2 \w_6}
 \oplustag{2} 2 V_{\w_2+\w_4}
 \oplustag{3} V_{\w_3+\w_5}
 \oplustag{3} 2 V_{\w_1+\w_2+\w_6}
 \oplustag{3} 3 V_{2 \w_2}
 \oplustag{4} V_{\w_4}
 \oplustag{2} V_{\w_4}
 \oplustag{3} V_{\w_1+\w_6}
 \oplustag{3} 2 V_{\w_2}
 \oplustag{3} V_{0}$
\item[$W_1(5)$] $\iso V_{\w_5}
 \oplustag{1} V_{\w_1}$
\item[$W_2(5)$] $\iso V_{2 \w_5}
 \oplustag{1} V_{\w_1+\w_5}
 \oplustag{2} V_{2 \w_1}$
\item[$W_3(5)$] $\iso V_{3 \w_5}
 \oplustag{1} V_{\w_1+2 \w_5}
 \oplustag{2} V_{2 \w_1+\w_5}
 \oplustag{3} V_{3 \w_1}$
\item[$W_m(6)$] remains irreducible for all $m$.

\item[{\boldmath $E_7$}]
\item[$W_1(1)$] $\iso V_{\w_1}
 \oplustag{1} V_{0}$
\item[$W_2(1)$] $\iso V_{2 \w_1}
 \oplustag{1} V_{\w_1}
 \oplustag{2} V_{0}$
\item[$W_1(2)$] $\iso V_{\w_2}
 \oplustag{1} V_{\w_7}$
\item[$W_2(2)$] $\iso V_{2 \w_2}
 \oplustag{1} V_{\w_2+\w_7}
 \oplustag{2} V_{2 \w_7}$
\item[$W_1(3)$] $\iso V_{\w_3}
 \oplustag{1} V_{\w_6}
 \oplustag{1} 2 V_{\w_1}
 \oplustag{1} V_{0}$
\item[$W_2(3)$] $\iso V_{2 \w_3}
 \oplustag{1} V_{\w_3+\w_6}
 \oplustag{2} V_{2 \w_6}
 \oplustag{1} 2 V_{\w_1+\w_3}
 \oplustag{2} V_{\w_4}
 \oplustag{2} 2 V_{\w_1+\w_6}
 \oplustag{2} 3 V_{2 \w_1}
 \oplustag{3} V_{\w_3}
 \oplustag{1} V_{\w_3}
 \oplustag{2} V_{\w_6}
 \oplustag{2} 2 V_{\w_1}
 \oplustag{2} V_{0}$
\item[$W_1(4)$] $\iso V_{\w_4}
 \oplustag{1} V_{\w_1+\w_6}
 \oplustag{1} 2 V_{\w_2+\w_7}
 \oplustag{1} V_{2 \w_1}
 \oplustag{1} 3 V_{\w_3}
 \oplustag{2} V_{\w_6}
 \oplustag{1} V_{2 \w_7}
 \oplustag{1} 3 V_{\w_6}
 \oplustag{2} V_{\w_1}
 \oplustag{1} 3 V_{\w_1}
 \oplustag{2} V_{0}
 \oplustag{1} V_{0}$
\item[$W_2(4)$] $\iso V_{2 \w_4}
 \oplustag{1} V_{\w_1+\w_4+\w_6}
 \oplustag{2} V_{2 \w_1+2 \w_6}
 \oplustag{1} 2 V_{\w_2+\w_4+\w_7}
 \oplustag{2} V_{\w_3+\w_5+\w_7}
 \oplustag{2} 2 V_{\w_1+\w_2+\w_6+\w_7}
 \oplustag{2} 3 V_{2 \w_2+2 \w_7}
 \oplustag{3} V_{\w_4+2 \w_7}
 \oplustag{1} V_{2 \w_1+\w_4}
 \oplustag{2} V_{3 \w_1+\w_6}
 \oplustag{2} V_{4 \w_1}
 \oplustag{1} 3 V_{\w_3+\w_4}
 \oplustag{2} 2 V_{\w_1+\w_2+\w_5}
 \oplustag{2} 4 V_{\w_1+\w_3+\w_6}
 \oplustag{2} V_{2 \w_5}
 \oplustag{2} V_{2 \w_2+\w_6}
 \oplustag{2} 3 V_{\w_4+\w_6}
 \oplustag{3} V_{\w_1+2 \w_6}
 \oplustag{2} 2 V_{2 \w_1+\w_2+\w_7}
 \oplustag{2} 6 V_{\w_2+\w_3+\w_7}
 \oplustag{3} 2 V_{\w_1+\w_5+\w_7}
 \oplustag{3} 2 V_{\w_2+\w_6+\w_7}
 \oplustag{2} 3 V_{2 \w_1+\w_3}
 \oplustag{2} 6 V_{2 \w_3}
 \oplustag{3} 3 V_{\w_1+\w_4}
 \oplustag{3} V_{\w_2+\w_5}
 \oplustag{3} V_{2 \w_1+\w_6}
 \oplustag{3} 3 V_{\w_3+\w_6}
 \oplustag{4} V_{2 \w_6}
 \oplustag{1} V_{\w_4+2 \w_7}
 \oplustag{2} V_{\w_1+\w_6+2 \w_7}
 \oplustag{2} 2 V_{\w_2+3 \w_7}
 \oplustag{2} V_{4 \w_7}
 \oplustag{1} 3 V_{\w_4+\w_6}
 \oplustag{2} 2 V_{\w_2+\w_3+\w_7}
 \oplustag{2} 3 V_{\w_1+2 \w_6}
 \oplustag{2} 4 V_{\w_1+\w_5+\w_7}
 \oplustag{2} V_{2 \w_3}
 \oplustag{2} V_{\w_1+2 \w_2}
 \oplustag{2} 3 V_{\w_1+\w_4}
 \oplustag{3} V_{2 \w_1+\w_6}
 \oplustag{2} 8 V_{\w_2+\w_6+\w_7}
 \oplustag{3} 2 V_{\w_3+2 \w_7}
 \oplustag{2} 6 V_{\w_2+\w_5}
 \oplustag{3} 2 V_{\w_3+\w_6}
 \oplustag{3} 2 V_{\w_1+\w_2+\w_7}
 \oplustag{2} V_{2 \w_1+2 \w_7}
 \oplustag{2} 3 V_{\w_3+2 \w_7}
 \oplustag{3} V_{\w_6+2 \w_7}
 \oplustag{2} 3 V_{2 \w_1+\w_6}
 \oplustag{3} V_{3 \w_1}
 \oplustag{2} 9 V_{\w_3+\w_6}
 \oplustag{3} 4 V_{\w_1+\w_2+\w_7}
 \oplustag{3} 3 V_{2 \w_6}
 \oplustag{3} 4 V_{\w_5+\w_7}
 \oplustag{3} 4 V_{\w_1+\w_3}
 \oplustag{3} V_{2 \w_2}
 \oplustag{3} 3 V_{\w_4}
 \oplustag{4} V_{\w_1+\w_6}
 \oplustag{2} 3 V_{\w_6+2 \w_7}
 \oplustag{2} 6 V_{2 \w_6}
 \oplustag{3} 3 V_{\w_5+\w_7}
 \oplustag{3} V_{\w_4}
 \oplustag{3} V_{\w_1+2 \w_7}
 \oplustag{3} 3 V_{\w_1+\w_6}
 \oplustag{4} V_{2 \w_1}
 \oplustag{1} 3 V_{\w_1+\w_4}
 \oplustag{2} 2 V_{\w_2+\w_5}
 \oplustag{2} 3 V_{2 \w_1+\w_6}
 \oplustag{2} 4 V_{\w_3+\w_6}
 \oplustag{2} 8 V_{\w_1+\w_2+\w_7}
 \oplustag{3} 2 V_{\w_5+\w_7}
 \oplustag{2} V_{2 \w_6}
 \oplustag{2} 3 V_{\w_5+\w_7}
 \oplustag{3} V_{\w_1+2 \w_7}
 \oplustag{2} 3 V_{3 \w_1}
 \oplustag{2} 12 V_{\w_1+\w_3}
 \oplustag{3} 4 V_{\w_4}
 \oplustag{3} 4 V_{\w_1+\w_6}
 \oplustag{2} 3 V_{2 \w_2}
 \oplustag{3} V_{\w_4}
 \oplustag{2} 5 V_{\w_4}
 \oplustag{3} 3 V_{\w_1+\w_6}
 \oplustag{3} 4 V_{\w_2+\w_7}
 \oplustag{3} V_{2 \w_1}
 \oplustag{3} 3 V_{\w_3}
 \oplustag{4} V_{\w_6}
 \oplustag{2} 3 V_{\w_1+2 \w_7}
 \oplustag{2} 9 V_{\w_1+\w_6}
 \oplustag{3} 4 V_{\w_2+\w_7}
 \oplustag{3} 3 V_{2 \w_1}
 \oplustag{3} 4 V_{\w_3}
 \oplustag{3} V_{2 \w_7}
 \oplustag{3} 3 V_{\w_6}
 \oplustag{4} V_{\w_1}
 \oplustag{2} 6 V_{2 \w_1}
 \oplustag{3} 3 V_{\w_3}
 \oplustag{3} V_{\w_6}
 \oplustag{3} 3 V_{\w_1}
 \oplustag{4} V_{0}
 \oplustag{1} V_{\w_4}
 \oplustag{2} V_{\w_1+\w_6}
 \oplustag{2} 2 V_{\w_2+\w_7}
 \oplustag{2} V_{2 \w_1}
 \oplustag{2} 3 V_{\w_3}
 \oplustag{3} V_{\w_6}
 \oplustag{2} V_{2 \w_7}
 \oplustag{2} 3 V_{\w_6}
 \oplustag{3} V_{\w_1}
 \oplustag{2} 3 V_{\w_1}
 \oplustag{3} V_{0}
 \oplustag{2} V_{0}$
\item[$W_1(5)$] $\iso V_{\w_5}
 \oplustag{1} V_{\w_1+\w_7}
 \oplustag{1} 2 V_{\w_2}
 \oplustag{1} 2 V_{\w_7}$
\item[$W_2(5)$] $\iso V_{2 \w_5}
 \oplustag{1} V_{\w_1+\w_5+\w_7}
 \oplustag{2} V_{2 \w_1+2 \w_7}
 \oplustag{1} 2 V_{\w_2+\w_5}
 \oplustag{2} V_{\w_3+\w_6}
 \oplustag{2} 2 V_{\w_1+\w_2+\w_7}
 \oplustag{2} 3 V_{2 \w_2}
 \oplustag{3} V_{\w_4}
 \oplustag{1} 2 V_{\w_5+\w_7}
 \oplustag{2} V_{\w_4}
 \oplustag{2} 2 V_{\w_1+2 \w_7}
 \oplustag{2} 2 V_{\w_1+\w_6}
 \oplustag{2} 4 V_{\w_2+\w_7}
 \oplustag{3} V_{\w_3}
 \oplustag{2} 3 V_{2 \w_7}
 \oplustag{3} V_{\w_6}$
\item[$W_1(6)$] $\iso V_{\w_6}
 \oplustag{1} V_{\w_1}
 \oplustag{1} V_{0}$
\item[$W_2(6)$] $\iso V_{2 \w_6}
 \oplustag{1} V_{\w_1+\w_6}
 \oplustag{2} V_{2 \w_1}
 \oplustag{1} V_{\w_6}
 \oplustag{2} V_{\w_1}
 \oplustag{2} V_{0}$
\item[$W_m(7)$] remains irreducible for all $m$.

\item[{\boldmath $E_8$}]
\item[$W_1(1)$] $\iso V_{\w_1}
 \oplustag{1} V_{\w_8}
 \oplustag{1} V_{0}$
\item[$W_1(2)$] $\iso V_{\w_2}
 \oplustag{1} V_{\w_7}
 \oplustag{1} 2 V_{\w_1}
 \oplustag{1} 2 V_{\w_8}
 \oplustag{1} V_{0}$
\item[$W_1(3)$] $\iso V_{\w_3}
 \oplustag{1} V_{\w_6}
 \oplustag{1} 2 V_{\w_1+\w_8}
 \oplustag{1} 3 V_{\w_2}
 \oplustag{2} V_{\w_7}
 \oplustag{1} V_{2 \w_8}
 \oplustag{1} 3 V_{\w_7}
 \oplustag{2} V_{\w_1}
 \oplustag{1} 4 V_{\w_1}
 \oplustag{2} 2 V_{\w_8}
 \oplustag{1} 3 V_{\w_8}
 \oplustag{2} V_{0}
 \oplustag{1} V_{0}$
\item[$W_1(4)$] $\iso V_{\w_4}
 \oplustag{1} V_{\w_1+\w_6}
 \oplustag{1} 2 V_{\w_2+\w_7}
 \oplustag{1} V_{2 \w_1+\w_8}
 \oplustag{1} 3 V_{\w_3+\w_8}
 \oplustag{2} V_{\w_6+\w_8}
 \oplustag{1} V_{2 \w_7}
 \oplustag{1} 6 V_{\w_1+\w_2}
 \oplustag{2} 2 V_{\w_5}
 \oplustag{2} 2 V_{\w_1+\w_7}
 \oplustag{1} 3 V_{\w_6+\w_8}
 \oplustag{2} V_{\w_1+2 \w_8}
 \oplustag{1} 5 V_{\w_5}
 \oplustag{2} 3 V_{\w_1+\w_7}
 \oplustag{2} 4 V_{\w_2+\w_8}
 \oplustag{2} V_{2 \w_1}
 \oplustag{2} 3 V_{\w_3}
 \oplustag{3} V_{\w_6}
 \oplustag{1} 3 V_{\w_1+2 \w_8}
 \oplustag{2} V_{3 \w_8}
 \oplustag{1} 9 V_{\w_1+\w_7}
 \oplustag{2} 4 V_{\w_2+\w_8}
 \oplustag{2} 3 V_{2 \w_1}
 \oplustag{2} 4 V_{\w_3}
 \oplustag{2} 4 V_{\w_7+\w_8}
 \oplustag{2} 3 V_{\w_6}
 \oplustag{3} V_{\w_1+\w_8}
 \oplustag{1} 10 V_{\w_2+\w_8}
 \oplustag{2} 4 V_{\w_3}
 \oplustag{2} 6 V_{\w_7+\w_8}
 \oplustag{2} 6 V_{\w_6}
 \oplustag{2} 8 V_{\w_1+\w_8}
 \oplustag{3} 2 V_{\w_2}
 \oplustag{1} 6 V_{2 \w_1}
 \oplustag{2} 3 V_{\w_3}
 \oplustag{2} V_{\w_6}
 \oplustag{2} 3 V_{\w_1+\w_8}
 \oplustag{3} V_{2 \w_8}
 \oplustag{1} 7 V_{\w_3}
 \oplustag{2} 5 V_{\w_6}
 \oplustag{2} 8 V_{\w_1+\w_8}
 \oplustag{2} 9 V_{\w_2}
 \oplustag{3} 3 V_{\w_7}
 \oplustag{2} V_{2 \w_8}
 \oplustag{2} 3 V_{\w_7}
 \oplustag{3} V_{\w_1}
 \oplustag{1} V_{3 \w_8}
 \oplustag{1} 8 V_{\w_7+\w_8}
 \oplustag{2} 3 V_{\w_6}
 \oplustag{2} 4 V_{\w_1+\w_8}
 \oplustag{2} 3 V_{2 \w_8}
 \oplustag{3} V_{\w_7}
 \oplustag{1} 7 V_{\w_6}
 \oplustag{2} 5 V_{\w_1+\w_8}
 \oplustag{2} 8 V_{\w_2}
 \oplustag{2} 3 V_{2 \w_8}
 \oplustag{2} 9 V_{\w_7}
 \oplustag{3} 3 V_{\w_1}
 \oplustag{2} 3 V_{\w_1}
 \oplustag{3} V_{\w_8}
 \oplustag{1} 15 V_{\w_1+\w_8}
 \oplustag{2} 8 V_{\w_2}
 \oplustag{2} 9 V_{2 \w_8}
 \oplustag{2} 12 V_{\w_7}
 \oplustag{2} 9 V_{\w_1}
 \oplustag{3} 3 V_{\w_8}
 \oplustag{2} 3 V_{\w_8}
 \oplustag{3} V_{0}
 \oplustag{1} 8 V_{\w_2}
 \oplustag{2} 6 V_{\w_7}
 \oplustag{2} 10 V_{\w_1}
 \oplustag{2} 6 V_{\w_8}
 \oplustag{1} 6 V_{2 \w_8}
 \oplustag{2} 3 V_{\w_7}
 \oplustag{2} V_{\w_1}
 \oplustag{2} 3 V_{\w_8}
 \oplustag{3} V_{0}
 \oplustag{1} 7 V_{\w_7}
 \oplustag{2} 5 V_{\w_1}
 \oplustag{2} 8 V_{\w_8}
 \oplustag{2} V_{0}
 \oplustag{1} 7 V_{\w_1}
 \oplustag{2} 5 V_{\w_8}
 \oplustag{2} 3 V_{0}
 \oplustag{1} 5 V_{\w_8}
 \oplustag{2} 3 V_{0}
 \oplustag{1} V_{0}$
\item[$W_1(5)$] $\iso V_{\w_5}
 \oplustag{1} V_{\w_1+\w_7}
 \oplustag{1} 2 V_{\w_2+\w_8}
 \oplustag{1} V_{2 \w_1}
 \oplustag{1} 3 V_{\w_3}
 \oplustag{2} V_{\w_6}
 \oplustag{1} 2 V_{\w_7+\w_8}
 \oplustag{1} 4 V_{\w_6}
 \oplustag{2} 2 V_{\w_1+\w_8}
 \oplustag{2} 2 V_{\w_2}
 \oplustag{1} 6 V_{\w_1+\w_8}
 \oplustag{2} 2 V_{\w_2}
 \oplustag{2} 2 V_{2 \w_8}
 \oplustag{2} 2 V_{\w_7}
 \oplustag{1} 5 V_{\w_2}
 \oplustag{2} 3 V_{\w_7}
 \oplustag{2} 4 V_{\w_1}
 \oplustag{1} 3 V_{2 \w_8}
 \oplustag{2} V_{\w_7}
 \oplustag{1} 5 V_{\w_7}
 \oplustag{2} 3 V_{\w_1}
 \oplustag{2} 4 V_{\w_8}
 \oplustag{1} 5 V_{\w_1}
 \oplustag{2} 3 V_{\w_8}
 \oplustag{2} V_{0}
 \oplustag{1} 4 V_{\w_8}
 \oplustag{2} 2 V_{0}
 \oplustag{1} V_{0}$
\item[$W_1(6)$] $\iso V_{\w_6}
 \oplustag{1} V_{\w_1+\w_8}
 \oplustag{1} 2 V_{\w_2}
 \oplustag{1} V_{2 \w_8}
 \oplustag{1} 3 V_{\w_7}
 \oplustag{2} V_{\w_1}
 \oplustag{1} 3 V_{\w_1}
 \oplustag{2} V_{\w_8}
 \oplustag{1} 3 V_{\w_8}
 \oplustag{2} V_{0}
 \oplustag{1} V_{0}$
\item[$W_1(7)$] $\iso V_{\w_7}
 \oplustag{1} V_{\w_1}
 \oplustag{1} 2 V_{\w_8}
 \oplustag{1} V_{0}$
\item[$W_1(8)$] $\iso V_{\w_8}
 \oplustag{1} V_{0}$

\fussy
\end{enumerate}

\chapter{Polynomial Relations Among Characters}
\label{chap_poly}

In this chapter we investigate some polynomial relations which appear to
hold among the characters of certain finite-dimensional representations
of $U_q(\ghat)$.  Conjecture~\ref{KRconj}, which was the basis for
Chapter~\ref{chap_decomp}, proposes a formula for these characters which
does indeed satisfy these relations.  (While this chapter is phrased in
the language of $U_q(\ghat)$, the statements in~\cite{KR} are all in the
language of Yangians; see the end of section~\ref{sec_algebras} for the
sad story of this translation.)

The main result of this chapter is that these polynomial relations have
only one solution, using a positivity condition on characters of
$U_q(\ghat)$.  Therefore a proof that the characters of $U_q(\ghat)$ do
indeed satisfy these relations would imply Conjecture~\ref{KRconj}.

\section{Introduction}
\label{sec_intro}

Retain the notions of the previous chapter: $\g$ is a complex
finite-dimensional simple Lie algebra, $\ghat$ its corresponding affine
Lie algebra.  Because of the inclusion of quantum enveloping algebras
$U_q(\g)\into U_q(\ghat)$, any finite-dimensional representation of
$U_q(\ghat)$ is a direct sum of irreducible representations of $U_q(\g)$.

Again we let $W_m(\l)$ denote the representation whose Drinfeld
polynomials were given in Definition~\ref{def_wml}, though this time it
is a representation of $U_q(\ghat)$.  Let $Q_m(\l)$ denote its character.
These characters appear to satisfy certain remarkable polynomial
identities.  When $\g$ is simply-laced, the identities have the simple
form
\begin{equation}
\label{polyrel-simple}
Q_m(\l)^2 = Q_{m-1}(\l)\,Q_{m+1}(\l) + \prod_{\l'\sim\l} Q_m(\l')
\end{equation}
for each $\l=1,\ldots,n$ and $m\geq1$.  The product is taken over all 
$\l'$ adjacent to $\l$ in the Dynkin diagram of $\g$.  Using these 
relations, it is possible to write any character $Q_m(\l)$ in 
terms of the characters $Q_1(\l)$ of the ``fundamental 
representations'' of $U_q(\ghat)$.

The main result of this chapter is that these equations have only one
solution where $Q_m(\l)$ is the character of a $U_q(\g)$-module with
highest weight $m\w_\l$; that is, where $Q_m(\l)$ is a positive integer
linear combination of irreducible $U_q(\g)$-characters with highest
weights sitting under $m\w_\l$.  (We restrict our attention to the
classical families $A_n$, $B_n$, $C_n$ and $D_n$.)  We use the polynomial
relations to write some of these multiplicities in terms of the
multiplicities in the characters $Q_1(\l)$, and the resulting
inequalities determine all of the multiplicities.

\section{Polynomial relations}
\label{sec_relations}

We will study the characters $Q_m(\l)$ of the finite-dimensional
representations $W_m(\l)$ of $U_q(\ghat)$, where $m=0,1,2,\ldots$ and
$\l=1,\ldots,n$, where $n=\rank(\g)$.  Since $U_q(\g)$ appears as a Hopf
subalgebra of $U_q(\ghat)$, we can talk about weights and characters of
$U_q(\ghat)$ modules by restricting our attention to the $U_q(\g)$
action.  From this point of view, $W_m(\l)$ has highest weight $m\w_\l$.
If $m=0$ then $W_m(\l)$ is the trivial representation and $Q_m(\l)=1$.
The objects $W_1(\l)$ and $Q_1(\l)$ are called the {\em fundamental}
representations and characters.

Let $V(\lambda)$ denote the character of the irreducible representation
of $U_q(\g)$ with highest weight $\lambda$.  We will write characters
$Q_m(\l)$ as sums $\sum m_\lambda V(\lambda)$.  We will refer to the
coefficients $m_\lambda$ as the multiplicity of $V(\lambda)$ in the sum.

The characters $Q_m(\l)$ when $\g$ is of type $A_n$ are known to satisfy
the so-called ``discrete Hirota relations.''  A conjectured
generalization of these relations appears in \cite{KR} for the classical
Lie algebras, and appear as the ``$Q$-system'' in \cite{Ku} for the
exceptional cases as well.  While we are only interested in the classical
cases, we will give the relations in full generality.

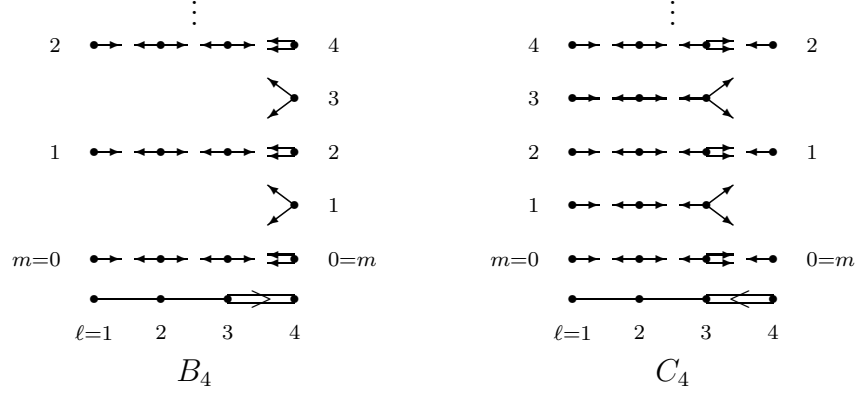
\begin{figure} \centering \setlength{\unitlength}{.035in}
\newcommand{\la}{{\vector(-1,0){4}}} \newcommand{\ra}{{\vector(1,0){4}}}
\newcommand{\ci}{{\circle*{1}}} \newcommand{\mpt}{\multiput}
\newsavebox{\lrbx}\sbox{\lrbx} {\begin{picture}(0,0)
  \mpt(0,0)(10,0){3}\ci\mpt(0,0)(10,0){2}\ra\mpt(10,0)(10,0){2}\la\end{picture}}
\newsavebox{\dd}\sbox{\dd} {\begin{picture}(0,0)
  \mpt(0,0)(10,0){4}\ci\put(0,0){\line(1,0){20}}
  \put(20,.6){\line(1,0){10}}\put(20,-.6){\line(1,0){10}}
\end{picture}}
\begin{picture}(70,60)
\mpt(20,20)(0,16){3}{\usebox{\lrbx}}
\mpt(40,20)(0,16){3}\ra
\mpt(50,20)(0,8){5}\ci
\mpt(50,20.6)(0,16){3}\la
\mpt(50,19.4)(0,16){3}\la
\mpt(50,28)(0,16){2}{\vector(-4,3){4}}
\mpt(50,28)(0,16){2}{\vector(-4,-3){4}} 
\put(20,14){\usebox{\dd}}
\put(45,14){\makebox(0,0){$>$}}
\put(20,9){\makebox(0,0){$\scriptstyle\l=1$}}
\put(30,9){\makebox(0,0){$\scriptstyle2$}}
\put(40,9){\makebox(0,0){$\scriptstyle3$}}
\put(50,9){\makebox(0,0){$\scriptstyle4$}}
\put(55,20){\makebox(0,0)[l]{$\scriptstyle 0=m$}}
\put(55,28){\makebox(0,0)[l]{$\scriptstyle 1$}}
\put(55,36){\makebox(0,0)[l]{$\scriptstyle 2$}}
\put(55,44){\makebox(0,0)[l]{$\scriptstyle 3$}}
\put(55,52){\makebox(0,0)[l]{$\scriptstyle 4$}}
\put(15,20){\makebox(0,0)[r]{$\scriptstyle m=0$}}
\put(15,36){\makebox(0,0)[r]{$\scriptstyle 1$}}
\put(15,52){\makebox(0,0)[r]{$\scriptstyle 2$}}
\put(35,3){\makebox(0,0){$B_4$}}
\put(35,58){\makebox(0,0){\vdots}}
\end{picture}
\begin{picture}(70,60)
\mpt(20,20)(0,8){5}{\usebox{\lrbx}}
\mpt(50,20)(0,16){3}\ci
\mpt(50,20)(0,16){3}\la
\mpt(40,20.6)(0,16){3}\ra
\mpt(40,19.4)(0,16){3}\ra
\mpt(40,28)(0,16){2}{\vector(4,3){4}}
\mpt(40,28)(0,16){2}{\vector(4,-3){4}} 
\put(20,14){\usebox{\dd}}
\put(45,14){\makebox(0,0){$<$}}
\put(20,9){\makebox(0,0){$\scriptstyle\l=1$}}
\put(30,9){\makebox(0,0){$\scriptstyle2$}}
\put(40,9){\makebox(0,0){$\scriptstyle3$}}
\put(50,9){\makebox(0,0){$\scriptstyle4$}}
\put(15,20){\makebox(0,0)[r]{$\scriptstyle m=0$}}
\put(15,28){\makebox(0,0)[r]{$\scriptstyle 1$}}
\put(15,36){\makebox(0,0)[r]{$\scriptstyle 2$}}
\put(15,44){\makebox(0,0)[r]{$\scriptstyle 3$}}
\put(15,52){\makebox(0,0)[r]{$\scriptstyle 4$}}
\put(55,20){\makebox(0,0)[l]{$\scriptstyle 0=m$}}
\put(55,36){\makebox(0,0)[l]{$\scriptstyle 1$}}
\put(55,52){\makebox(0,0)[l]{$\scriptstyle 2$}}
\put(35,3){\makebox(0,0){$C_4$}}
\put(35,58){\makebox(0,0){\vdots}}
\end{picture}
\caption{Product from the definition of $\cal Q$ for $B_4$ and $C_4$}
\label{fig_product}
\end{figure}

For every positive integer $m$ and for $\l=1,\ldots,n$,
\begin{equation}
\label{polyrel}
Q_m(\l)^2 = Q_{m+1}(\l)\,Q_{m-1}(\l) + \prod_{\l'\sim\l} {\cal Q}(m,\l,\l')
\end{equation}
The product is over all $\l'$ adjacent to $\l$ in the Dynkin diagram of
$\g$, and the contribution ${\cal Q}(m,\l,\l')$ from $\l'$ is determined
by the relative lengths of the roots $\alpha_\l$ and $\alpha_{\l'}$, as
follows:
\begin{equation}
\label{Qcurly}
\newcommand{\lal}{{\langle \alpha_\l,\alpha_\l \rangle}}
\newcommand{\lalp}{{\langle \alpha_{\l'},\alpha_{\l'} \rangle}}
\newlength{\kw} \settowidth{\kw}{$k$} \newcommand{\nok}{{\hspace{\kw}}}
{\cal Q}(m,\l,\l') =
 \left\{ \begin{array}{ll}
\displaystyle Q_m(\l') &
  \mbox{if \,} \nok\lal = \nok\lalp \\ [5pt]
\displaystyle Q_{km}(\l') &
  \mbox{if \,} \nok\lal = k\lalp \\ [2pt]
\displaystyle \prod_{i=0}^{k-1} Q_{\lfloor \frac{m+i}{k} \rfloor}(\l') &
  \mbox{if \,} k\lal = \nok\lalp
\end{array} \right.
\end{equation}
where $\lfloor x\rfloor$ is the greatest integer not exceeding $x$.  We
note that in the classical cases, the product differs from the simplified
version in equation~(\ref{polyrel-simple}) only when:
$$
\begin{array}{rll}
\g=\mathfrak{so}(2n+1),& \l=n-1: & Q_m(n-2)\,Q_{2m}(n) \\
&\l=n: & Q_{\lfloor\frac{m}{2}\rfloor}(n-1)\,
        Q_{\lfloor\frac{m+1}{2}\rfloor}(n-1) \\
\g=\mathfrak{sp}(2n),& \l=n-1: &  Q_m(n-2)\,
        Q_{\lfloor\frac{m}{2}\rfloor}(n)\,
        Q_{\lfloor\frac{m+1}{2}\rfloor}(n) \\
&\l=n: & Q_{2m}(n-1)
\end{array}
$$
The structure of the product is easily represented graphically, with a 
vertex for each character $Q_m(\l)$ and an arrow from $Q_m(\l)$ 
pointing at each term of $\prod {\cal Q}(m,\l,\l')$; see 
Figure~\ref{fig_product} for $\g$ of type $B_4$ and $C_4$.  The 
corresponding picture for $G_2$ is similarly pleasing.

Finally, we can solve equation~(\ref{polyrel}) to get a recurrence relation:
\begin{equation}
\label{recurrence}
Q_{m}(\l) = 
\frac{ Q_{m-1}(\l)^2 - \prod {\cal Q}(m-1,\l,\l') } {Q_{m-2}(\l)}
\end{equation}

Note that the recurrence is well-founded: repeated use eventually writes
everything in terms of the fundamental characters $Q_1(\l)$.  This is
just the statement that iteration of ``move down, then follow any arrow''
in Figure~\ref{fig_product} will eventually lead you from any point to
one on the bottom row.  In fact, $Q_m(\l)$ is always a polynomial in the
fundamental characters, though from looking at the recurrence it is only
clear that it is a rational function.  A Jacobi-Trudi style formula for
writing the polynomial directly was given in \cite{Hi}.

The reason that characters of representations of quantum affine algebras
are solutions to a discrete integrable system is still a bit of a
mystery.

\section{Main Theorem}
\label{sec_statement}

The result of \cite{KR} was to conjecture a combinatorial formula for all
the multiplicities $Z(m,\l,\lambda)$ in the decomposition $Q_m(\l)=\sum
Z(m,\l,\lambda) V(\lambda)$, reproduced as Conjecture~\ref{KRconj} here.
We will refer to these proposed characters as ``combinatorial
characters'' of the representations $W_m(\l)$, although it is unproven
that they are the characters of some $U_q(\ghat)$ module.

\begin{thm}[Kirillov-Reshetikhin]
\label{thm_KR}
Let $\g$ be of type $A$, $B$, $C$ or $D$.  The combinatorial 
characters of $W_m(\l)$ are the unique solution to 
equations~(\ref{polyrel}) and~(\ref{Qcurly}) with the initial data
$$
\begin{array}{llcll}
A_n:
& Q_1(\l) &=& V(\w_\l) & 1\leq\l\leq n \\
B_n:
& Q_1(\l) &=& V(\w_\l)+V(\w_{\l-2})+V(\w_{\l-4})+\cdots & 1\leq\l\leq n-1 \\
& Q_1(n)  &=& V(\w_n) \\
C_n:
& Q_1(\l) &=& V(\w_\l) & 1\leq\l\leq n \\
D_n:
& Q_1(\l) &=& V(\w_\l)+V(\w_{\l-2})+V(\w_{\l-4})+\cdots & 1\leq\l\leq n-2 \\
& Q_1(\l) &=& V(\w_\l) & \l=n-1,n
\end{array}
$$
\end{thm}
The main result of this chapter is that the specification of initial 
data is unnecessary.

\begin{thm}
\label{thm_main}
Let $\g$ be of type $A$, $B$, $C$ or $D$.  The combinatorial 
characters of the representations $W_m(\l)$ are the only solutions to 
equations~(\ref{polyrel}) and~(\ref{Qcurly}) such that $Q_m(\l)$ is a 
character of a representation of $U_q(\g)$ with highest weight 
$m\w_\l$, for every nonnegative integer $m$ and $1\leq\l\leq n$.
\end{thm}

We need only prove that any choice of initial data other than that in
Theorem~\ref{thm_KR} would result in some $Q_m(\l)$ which is not a
character of a representation of $U_q(\g)$.  The values $Q_m(\l)$ are
always virtual $U_q(\g)$-characters, but in all other cases, some contain
representations occurring with negative multiplicity.  As an immediate
consequence, we have:

\begin{cor}
\label{cor_main}
If the characters $Q_m(\l)$ of the representations $W_m(\l)$ obey the 
recurrence relations in equations~(\ref{polyrel}) and~(\ref{Qcurly}), then 
they must be given by the formula for combinatorial characters in 
Conjecture~\ref{KRconj}.
\end{cor}

The technique for proving Theorem~\ref{thm_main} is as follows.  The 
possible choices of initial data are limited by the requirement that 
$Q_1(\l)$ be a representation with highest weight $\w_\l$.  That is, 
$Q_1(\l)$ must decompose into irreducible $U_q(\g)$-modules as
$$
Q_1(\l) = V(\w_\l) + \sum_{\lambda\prec\w_\l} m_\lambda V(\lambda)
$$
Note that we require that $V(\w_\l)$ occur in $Q_1(\l)$ exactly once.  
Furthermore, we require that for every other component $V(\lambda)$ 
that appears, $\lambda\prec\w_\l$, {\em i.e.} that $\w_\l-\lambda$ is 
a positive root.

We proceed with a case-by-case proof.  For each series, we find 
explicit multiplicities of irreducible representations occurring in 
$Q_m(\l)$ which would be negative for any choice of $Q_1(\l)$ other 
than that of Theorem~\ref{thm_KR}.  The calculations for series $B$, 
$C$ and $D$ are found below.

When $\g$ is of type $A_n$, no computations are necessary, because 
every fundamental root is minuscule: there are no $\lambda\prec\w_\l$ 
to worry about, no other choices for initial data to rule out.  In 
fact, $Q_m(\l)$ is just $V(m\w_\l)$ for all $m$ and $\l$, and moreover 
every $U_q(\g)$ module is also acted upon by $U_q(\ghat)$, by means of 
the evaluation representation.

\subsection*{\boldmath $B_n$}

Let $\g$ be of type $B_n$.  Let $V_i$ stand for $V(\w_i)$ for $1\leq
i\leq n-1$, and $V_{sp}$ for the character of the spin representation
with highest weight $\w_n$.  For convenience, let $\w_0=0$ and $V_0$
denote the character of the trivial representation.  Finally, we denote
by $V_n$ the character of the representation with highest weight $2\w_n$,
which behaves like the fundamental representations.

There are no dominant weights $\lambda\prec\w_n$, so $Q_1(n)=V_{sp}$.
The only weights $\lambda\prec\w_a$ are $0,\w_1,\ldots,\w_{a-1}$ for
$1\leq a\leq n-1$, so we write
\begin{equation}
Q_1(a) = V_a + \sum_{b=0}^{a-1} M_{a,b} V_b
\end{equation}
Our goal is to prove that the only possible values for the
multiplicities are
\begin{equation}
\label{BMs}
M_{a,b} = \left\{
\begin{array}{ll} 1, & a-b \mbox{ even} \\ 0, & a-b \mbox{ odd} \end{array}
\right.
\end{equation}

We will show these values are necessary inductively; the proof for
each $M_{a,b}$ will assume the result for all $M_{c,d}$ with
$\lceil\frac{c-d}2\rceil < \lceil\frac{a-b}2\rceil$ as well as those
with $\lceil\frac{c-d}2\rceil = \lceil\frac{a-b}2\rceil$ and
$c+d>a+b$.  (Here $\lceil x \rceil$ is the least integer greater than or
equal to $x$.) This amounts to working in the following order:
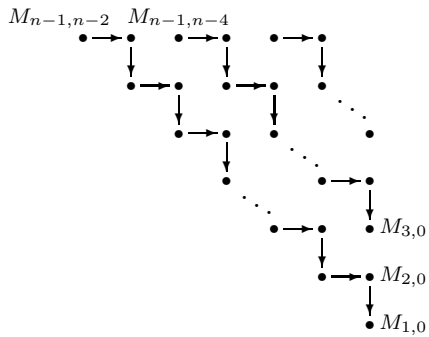
\begin{center} \setlength{\unitlength}{.25in}
\newcommand{\vx}{{\circle*{.15}}}
\newcommand{\rt}{{\vector(1,0){.6}}}
\newcommand{\dn}{{\vector(0,-1){.6}}}
\begin{picture}(8,7)
\multiput(0,6)(1,-1){7}{\vx}
\multiput(1,6)(1,-1){3}{\vx}
\multiput(5,2)(1,-1){2}{\vx}
\multiput(2,6)(1,-1){5}{\vx}
\multiput(3,6)(1,-1){2}{\vx}
\put(6,3){\vx}
\multiput(4,6)(1,-1){3}{\vx}
\put(5,6){\vx}
\multiput(0.2,6)(1,-1){3}{\rt}
\multiput(1,5.8)(1,-1){3}{\dn}
\multiput(4.2,2)(1,-1){2}{\rt}
\multiput(5,1.8)(1,-1){2}{\dn}
\multiput(2.2,6)(1,-1){2}{\rt}
\multiput(3,5.8)(1,-1){2}{\dn}
\put(5.2,3){\rt} \put(6,2.8){\dn}
\put(4.2,6){\rt} \put(5,5.8){\dn}
\multiput(3.25,2.25)(1,1){3}{\mbox{$\ddots$}}
\put(-.5,6.2){\makebox(0,0)[b]{$\scriptstyle M_{n-1,n-2}$}}
\put(2,6.2){\makebox(0,0)[b]{$\scriptstyle M_{n-1,n-4}$}}
\put(6.2,0){\makebox(0,0)[l]{$\scriptstyle M_{1,0}$}}
\put(6.2,1){\makebox(0,0)[l]{$\scriptstyle M_{2,0}$}}
\put(6.2,2){\makebox(0,0)[l]{$\scriptstyle M_{3,0}$}}
\end{picture}
\begin{minipage}[b]{2in}\raggedright
First follow the diagonal from $M_{n-1,n-2}$ to $M_{1,0}$, then the one
from $M_{n-1,n-4}$ to $M_{3,0}$, etc., ending in the top right corner
with $M_{n-1,0}$ or $M_{n-2,0}$, depending on the parity of $n$.
\vspace{12pt}
\end{minipage}
\end{center}
We show that equation~(\ref{BMs}) must hold for $M_{a,b}$, assuming it 
holds for all earlier $M$s in this ordering, by the following 
calculations:
\begin{enumerate}
\item
For $M_{n-1,b}$ where $n-1-b$ is odd, the multiplicity of $V(\w_b+\w_n)$
in $Q_3(n)$ is $1-2M_{n-1,b}$,
\item
For $M_{a,b}$ where $a-b$ is odd and $a\leq n-2$, the multiplicity of
$V(\w_{a+2}+\w_b)$ in $Q_2(a+1)$ is $-M_{a,b}$,
\item
For $M_{a,b}$ where $a-b$ is even:
\begin{itemize}
\item The multiplicity of $V(\w_a+\w_b)$ in $Q_2(a)$ is $2M_{a,b}-1$, and
\item The multiplicity of $V(\w_{a+2}+\w_b)$ in $Q_2(a+1)$ is $1-M_{a,b}$.
\end{itemize}
\end{enumerate}
Since all $M_{a,b}$ and all multiplicities are nonnegative integers, we
must have $M_{a,b}=0$ to satisfy the first two cases and $M_{a,b}=1$ to
satisfy the third.

The calculations to prove these claims depend on the ability to tensor
together the $U_q(\g)$-modules whose characters form $Q_m(\l)$.  A
complete algorithm for decomposing these tensors is given in terms of
crystal bases in \cite{N}.  For the current case, though, it happens that
the only tensors we need to take are of fundamental representations.
Simple explicit formulas for these decompositions had been given in
\cite{KN} before the advent of crystal base technology.

\begin{enumerate}
\item $M_{n-1,b}$, $n-1-b$ odd:

We want to find the multiplicity of $V(\w_b+\w_n)$ in $Q_3(n)$.
Recursing through the polynomial relations, we find that
$$
Q_3(n) = Q_1(n)^3 - 2 Q_1(n) Q_1(n-1) =
Q_1(n) \left[ Q_1(n)^2-2Q_1(n-1) \right]
$$
Assuming equation~(\ref{BMs}) for $M_{n-1,b'}$ for $b'>b$ and recalling that
$Q_1(n)^2 = V_{sp}^2 = V_n+V_{n-1}+\cdots+V_0$, we need to compute
the product
$$
V_{sp}  \left[
V_n-V_{n-1}+V_{n-2}-\cdots-V_{b+1}+(1-2M_{n-1,b})V_b-\cdots \right]
$$
Since $V_{sp}V_k=\sum_{i=0}^kV(\w_i+\w_n)$, we find that the
multiplicity of $V(\w_b+\w_n)$ in the product is the desired
$1-2M_{n-1,b}$.

\item $M_{a,b}$, $a-b$ odd, $a\leq n-2$:

This calculation is typical of many of the ones that will follow, and
will be written out in more detail.  We want to know the multiplicity of
$V(\w_{a+2}+\w_b)$ in $Q_2(a+1)$.  When $a\leq n-3$, we have
$$
Q_2(a+1)=Q_1(a+1)^2 - Q_1(a+2) Q_1(a)
$$
Assuming equation~(\ref{BMs}) holds for all earlier $M$s in the 
ordering, we have
\begin{eqnarray*}
Q_1(a+1) &=& V_{a+1} + V_{a-1} + \cdots+V_{b}+M_{a+1,b-1}V_{b-1}+\cdots \\
Q_1(a+2) &=& V_{a+2} + V_{a} + \cdots+V_{b+1}+M_{a+2,b}V_{b}+\cdots \\
Q_1(a)   &=& V_{a} + V_{a-2} + \cdots+V_{b+1}+M_{a,b}V_{b}+\cdots
\end{eqnarray*}
To compute $Q_1(a+1)^2 - Q_1(a+2) Q_1(a)$, we note that the
$V_s V_t$ term in $Q_1(a+1)^2$ and the $V_{s+1} V_{t-1}$
term in $Q_1(a+2) Q_1(a)$ are almost identical: when $s>t$, for example,
the difference is just $\sum_{i=0}^t V(\w_i+\w_{s-t-2+i})$.  In our
case, the only $V(\w_{a+2}+\w_b)$ term that does not cancel out is the
one contributed by $M_{a,b} V_{a+2} V_b$, and the
multiplicity of $V(\w_{a+2}+\w_b)$ is $-M_{a,b}$.

When $a=n-2$ the polynomial relations instead look like
$$
Q_2(n-1)=Q_1(n-1)^2 - Q_1(n)^2 Q_1(n-2) + Q_1(n-1) Q_1(n-2)
$$
The $V_n+V_{n-2}+\cdots$ terms of $Q_1(n)^2$ behave just like the
$Q_1(a+2)$ term in the above argument.  The extra terms from
$Q_1(n-2)\left[Q_1(n-1)-V_{n-1}-V_{n-3}-\cdots\right]$ make no net
contribution, as can be seen by checking highest weights.

\item $M_{a,b}$, $a-b$ even:

Calculating the multiplicity of $V(\w_{a+2}+\w_b)$ in $Q_2(a+1)$ is
similar to the above; the trick of canceling $V_s V_t$
with $V_{s+1} V_{t-1}$ works again.  The only terms remaining are
$+1$ from $V_{a+1} V_{b+1}$ and the same $-M_{a,b}$ from
$V_{a+2} M_{a,b}V_b$ as above, so the multiplicity is $1-M_{a,b}$

Likewise, calculating the multiplicity of $V(\w_a+\w_b)$ in $Q_2(a)$ we
find two contributions of $M_{a,b}$ from $M_{a,b}V_a V_b$ (in either
order) in $Q_2(a)^2$, and a contribution of $1$ from $V_{a+1} V_{b+1}$ in
$Q_1(a+1) Q_1(a-1)$, so the multiplicity is $2M_{a,b}-1$.
\end{enumerate}

\subsection*{\boldmath $C_n$}

Let $\g$ be of type $C_n$.  We let $V_i$ stand for $V(\w_i)$ for $1\leq
i\leq n$.  The only dominant weights $\lambda\prec\w_a$ for $1\leq a\leq
n$ are $\lambda=\w_b$ for $0\leq b<a$ and $a-b$ even, where $\w_0=0$.
(If $a-b$ is odd, then $\w_a$ and $\w_b$ lie in different translates of
the root lattice, so are incomparable.)  So we write
\begin{equation}
Q_1(a) = V_a + \sum_{i=0}^{\lfloor a/2 \rfloor} M_{a,a-2i} V_{a-2i}
\end{equation}
We will prove that in fact $M_{a,b}=0$ for all $a$ and $b$.

Again we choose a convenient order to investigate the multiplicities:
first look at $M_{a,a-2}$ for $a=n,n-1,\ldots,2$, and then all $M_{a,b}$
with $a-b=4,6,8,\ldots$.  This time the multiplicities acting as
witnesses are:
\begin{enumerate}
\item
For $M_{a,a-2}$, the multiplicity of $V(\w_{a-1}+2\w_{a-2})$ in
$Q_3(a-1)$ is $1-2M_{a,a-2}$,
\item
For $M_{a,b}$ for $a-b\geq4$, the multiplicity of $V(\w_{a-2}+\w_b)$
in $Q_2(a-1)$ is $-M_{a,b}$.
\end{enumerate}
Performing these computations requires the ability to tensor more general
representations of $\g$ than were needed in the $B_n$ case.  For this we
use the generalization of the Littlewood-Richardson rule to all classical
Lie algebras given in \cite{N}, which we summarize briefly in an Appendix
to this chapter.

\pagebreak[2]
\begin{enumerate}
\item $M_{a,a-2}$:

We want to calculate the multiplicity of $V(\w_{a-1}+2\w_{a-2})$ in
$Q_3(a-1)$.  First we write $Q_3(a-1)$ as a sum of terms of the form
$Q_1(x)Q_1(y)Q_1(z)$, which we denote as $(x;y;z)$ for brevity.  When
$2\leq a-1\leq n-2$, we have
\begin{eqnarray*}
Q_3(a-1) &=&
(a-1;a-1;a-1)-2(a;a-1;a-2)\\
&&{}-(a+1;a-1;a-3)+(a;a;a-3)+(a+1;a-2;a-2)
\end{eqnarray*}

When $a-1$ is one of $1,2$ or $n-1$, the above decomposition still 
holds, if we set $Q_1(0)=1$ and $Q_1(-1)=Q_1(n+1)=0$.  We want to find 
the multiplicity of $V(\w_{a-1}+2\w_{a-2})$ in each of these terms.

First, $V(\w_{a-1}+2\w_{a-2})$ occurs with multiplicity 3 in the
$V_{a-1}^3$ component of $Q_1(a-1)^3$.  We calculate this
number using the crystal basis technique for tensoring representations.
Beginning with the Young diagram of $V_{a-1}$, we must choose a tableau
$1,2,\ldots,a-2,p$ from the second tensor factor, where $p$ must be be
one of $a-1$, $a$, or $\overline{a-1}$.  Then the choice of tableau from
the third tensor component must be the same but replacing $p$ with
$\overline{p}$.

Similarly, the $V_{a} V_{a-1} V_{a-2}$ component of the
$(a;a-1;a-2)$ term produces $V(\w_{a-1}+2\w_{a-2})$ with multiplicity 1,
corresponding to the choice of the tableau $1,2,\ldots,a-2,\overline{a}$
from the crystal of $V_{a-1}$.  We see that the remaining three terms
cannot contribute by looking at tableaux in the same way.

Second, $V(\w_{a-1}+2\w_{a-2})$ occurs in the $M_{a,a-2}V_{a-2} 
V_{a-1} V_{a-2}$ piece of $(a;a-1;a-2)$ and the 
$M_{a+1,a-1}V_{a-1} V_{a-2} V_{a-2}$ piece of 
$(a+1;a-2;a-2)$ as the highest weight component.  Our inductive 
hypothesis, however, assumes that $M_{a+1,a-1}=0$, and we start the 
induction with $a=n$, where the $(a+1;a-2;a-2)$ term vanishes 
entirely.

Totaling these results, we find that the net multiplicity is 
$1-2M_{a,a-2}$, and conclude that $M_{a,a-2}=0$.

\item $M_{a,b}$ for $a-b\geq4$:

We want to calculate the multiplicity of $V(\w_{a-2}+\w_b)$
in $Q_2(a-1)$.  For any $2\leq a-1\leq n-1$, we have
\begin{eqnarray*}
Q_2(a-1) &=& Q_1(a-1)^2 - Q_1(a)Q_1(a-2) \\
&=& (V_{a-1}+\cdots)(V_{a-1}+\cdots) 
  - (V_a+M_{a,b}V_b+\cdots)(V_{a-2}+\cdots)
\end{eqnarray*}
where every omitted term is either already known to be 0 by induction, or
else has highest weight less than $\w_b$, so cannot contribute.  As in
the $B_n$ case, the $V_{a-1}^2$ and $V_aV_{a-2}$ terms nearly cancel one
another's contributions: their difference is just $\sum_{k=0}^{a-1}V(2\w_k)$.
Since $V(\w_{a-2}+\w_b)$ occurs in $V_{a-2}V_b$ with multiplicity 1, the
net multiplicity in $Q_2(a-1)$ is $-M_{a,b}$, and we conclude that
$M_{a,b}=0$.
\end{enumerate}

\subsection*{\boldmath $D_n$}

Let $\g$ be of type $D_n$.  This time we let $V_i$ stand for $V(\w_i)$
for $1\leq i\leq n-2$, and use $V_{n-1}$ for the character of the
representation with highest weight $\w_{n-1}+\w_n$.  We will not need to
explicitly use the characters of the two spin representations
individually, only their product, $V_{n-1}+V_{n-3}+\cdots$.

There are no dominant weights under $\w_{n-1}$ or $\w_n$, and so no work
to do on $Q_1(n-1)$ or $Q_1(n)$.  For $1\leq a\leq n-2$, the only
dominant weights $\lambda\prec\w_a$ are $\lambda=\w_b$ for $0\leq b<a$
and $a-b$ even; again $\w_0=0$.  (If $a-b$ is odd, then $\w_a$ and $\w_b$
lie in different translates of the root lattice, so are incomparable.)
So we write
\begin{equation}
Q_1(a) = V_a + \sum_{i=0}^{\lfloor a/2 \rfloor} M_{a,a-2i} V_{a-2i}
\end{equation}
We will show that in fact $M_{a,b}=1$ for all $a$ and $b$.

Again the proof is by induction; to show $M_{a,b}=1$ we will assume
$M_{c,d}=1$ as long as either $c-d<a-b$ or $c-d=a-b$ and $c>a$.  (This is
the same ordering used for the $B_n$ series after dropping the $M_{a,b}$
with $a-b$ odd.)  Our witnesses this time are:
\begin{itemize}
\item The multiplicity of $V(2\w_b)$ in $Q_2(a-1)$ is $1-M_{a,b}$, and
\item The multiplicity of $V(\w_a+\w_b)$ in $Q_2(a)$ is $2M_{a,b}-1$.
\end{itemize}
We must therefore conclude that $M_{a,b}=1$.  Since we only need to
tensor fundamental representations together, the explicit formulas given
in \cite{KN} are enough to carry out these calculations.

For any $\l\leq n-3$, the polynomial relations give us
$$
Q_2(\l)=Q_1(\l)^2-Q_1(\l+1)Q_1(\l-1)
$$

The multiplicity of $V(2\w_b)$ in $Q_2(a-1)$ is easily calculated
directly, since $V(2\w_b)$ appears in $V_r V_s$ if and only if
$r=s\geq b$, and then it appears with multiplicity one.  The $Q_1(a-1)$
term therefore contains $V(2\w_b)$ exactly $(a-b)/2$ times, while the
$Q_1(a)Q_1(a-2)$ term subtracts off $M_{a,b}-1+(a-b)/2$ of them.
Thus the net multiplicity is $1-M_{a,b}$.

To calculate the multiplicity of $V(\w_a+\w_b)$ in $Q_2(a)$ for $a\leq
n-3$, we once again use the trick of canceling the contribution from the
$V_s V_t$ term of $Q_1(a)^2$ with the $V_{s+1} V_{t-1}$ term of
$Q_1(a+1)Q_1(a-1)$.  The cancellation requires more attention this time,
since $V(\w_a+\w_b)$ occurs with multiplicity two in $V_s V_t$ when
$a-b\geq 2n-r-s$.  In the end, the only terms that do not cancel are the
contributions of $M_{a,b}$ from $V_a V_b$ and $V_b V_a$ in $Q_1(a)^2$ and
of $-1$ from $V_{b+1} V_{a-1}$ in $Q_1(a+1)Q_1(a-1)$.  Thus the net
multiplicity is $2M_{a,b}-1$.

Finally, if $a=n-2$ the polynomial relations change to
$$
Q_2(n-2)=Q_1(n-2)^2-Q_1(n-1)Q_1(n)Q_1(n-3)
$$
This change does not require any new work, though: $Q_1(n-1)Q_1(n)$ is
just the product of the two spin representations, which decomposes as
$V_{n-1}+V_{n-3}+\cdots$.  Since this is exactly what we wanted
$Q_1(\l+1)$ to look like in the above argument, the preceding calculation
still holds.

\section*{Appendix: Littlewood-Richardson Rule for $C_n$}

This is a brief summary of a generalization of the Littlewood-Richardson
rule to Lie algebras of type $C_n$, as given in \cite{N}.  For our
purposes, we only need the ability to tensor an arbitrary representation
with one of the fundamental representations with highest weights
$\w_1,\ldots,\w_n$.

The representation with highest weight $\sum_{k=1}^{n} a_k\w_k$ is
represented by a Young diagram $Y$ with $a_k$ columns of height $k$.  For
a fundamental representation $V_k$, we create Young tableaux from our
column of height $k$ by filling in the boxes with $k$ distinct symbols
$i_1,\ldots,i_k$ chosen in order from the sequence
$1,2,\ldots,n,\overline{n},\ldots,\overline{2},\overline{1}$ in all
possible ways, as long as if $i_a=p$ and $i_b=\overline{p}$ then
$a+(k-b+1)\leq p$.  These tableaux label the vertices of the crystal
graph of the representation $V_k$.

Given a Young diagram $Y$, the symbols $1,2,\ldots,n$ act on it by adding
one box to the first, second,\ldots,$n$th row, and the symbols
$\overline{1},\overline{2},\ldots,\overline{n}$ act by removing one,
provided the addition or removal results in a diagram whose rows are
still nonincreasing in length.  The result of the action of the symbol
$i_a$ on $Y$ is denoted $Y\!\ba i_a$.

Then the tensor product $V\tensor V_k$, where $V$ has Young diagram $Y$,
decomposes as the sum of all representations with diagrams
$(((Y\!\ba i_1)\ba i_2)\cdots\ba i_k)$, where $i_1,\ldots,i_k$ range over
all tableaux of $V_k$ such that each of the actions still result in a
diagram whose rows are still nonincreasing in length.

\chapter{Summary and Further Questions}
\label{chap_questions}

Our goal was to study finite-dimensional representations of Yangians and
quantum affine algebras.  We began with a conjectural formula of Kirillov
and Reshetikhin's, based on methods of mathematical physics, describing
how some of these representations should decompose into irreducible
representations of the underlying finite-dimensional Lie algebras.

In Chapter~\ref{chap_decomp}, we gave a new combinatorial interpretation
of this formula.  This new point of view made possible computations which
were completely intractable using the original version.  This formulation
also endows the decomposition into irreducibles with a new tree
structure, if the underlying Lie algebra is simply-laced.  We used this
tree to note some structural compatibility among representations whose
highest weights were different multiples of the same fundamental weight,
and were able to calculate the asymptotics of the growth of their
dimension.

In Chapter~\ref{chap_poly}, we explored a set of polynomial relations
that seem to hold among the characters of some of these representations.
The fact that every finite-dimensional representation of the quantum
affine algebra is a direct sum of representations of the underlying Lie
algebra is a very strong positivity condition.  We proved that for the
classical families of Lie algebras, the positivity condition and the
polynomial relations leave only one choice for the characters of the
quantum affine algebras.

We conclude with some questions for further research which seem to come
naturally from the topics discussed here.

\begin{enumerate}

\item
First and foremost, the conjectural Kirillov-Reshetikhin formula requires
a solid mathematical proof.  In light of our results in
Chapter~\ref{chap_poly}, it would suffice to prove that the characters of
the representations $\wml$ (as $\g$-modules) actually do satisfy the
polynomial relations.  One possible way to prove this would be by showing
the representations form an exact sequence, something like
$$
0 \to
W_{m+1}(\l)\tensor W_{m-1}(\l) \to
W_{m}(\l)\tensor W_{m}(\l) \to
\bigotimes_{\l'\sim\l} {\cal W}(m,\l,\l') \to
0
$$
If we want to treat these as tensor products of $U_q(\ghat)$-modules, we
need to specify shifts as well.  One can hope that the correct shifts are
given by the ``T-system'' written down in~\cite{Ku}, a proposed
generalization of these polynomial relations to include shifts.  It is
also possible that there are two different choices for the shifts such
that the above exact sequence is correct for one and backwards for the
other.

\item
It seems reasonable to hope that these results apply to cases other than
just $U_q(\ghat)$.  Certainly the twisted simply-laced affine cases are
the most natural candidates.  One can also hope that the polynomial
relations hold if we replace $\g$ with a more general Kac-Moody algebra.
Of course, the easy definition of $U_q(\ghat)$ is no longer available.
But since the results of Chapter~\ref{chap_poly} do not rely on knowing
what representations of $U_q(\ghat)$ look like, it may be possible to
generalize the polynomial relations and find a unique solution first, and
only later identify the solutions as irreducible representations of some
new object.

\item
The tree structure introduced here is, so far, only a little better than
a computational tool.  The fact that it highlighted a similarity in
structure between representations with different highest weights, though,
indicates that there may be some representation theoretic meaning to the
way it arranges pieces in the decomposition.  If so, there might be a
similar structure to decompositions for non-simply-laced cases, even
though this construction does not make sense.

A probable first step would be to understand better the structure in the
$A_n$ case, where a body of knowledge exists on tensor products of
rectangles.  Discretion is required, though: any patterns noticed for
$A_n$ can generalize either to the other classical finite-dimensional Lie
algebras or to representations of Yangians.

\item
The original Kirillov-Reshetikhin formula, or even just the tree
structure in the $A_n$ case, should have a generalization to
representations with highest weights other than just multiples of a
fundamental weight.  Unfortunately, the algorithm for tensoring
rectangles, which is phrased in a way that would make sense for tensoring
arbitrary shapes, is certainly not true beyond the rectangle case: as
presented here, it only depends on the multisets of $n$th columns of its
arguments, which the tensor product does not.

There is another potential stumbling block in generalizing these results
in cases other than $A_n$: it is no longer clear which 
$U_q(\g)$-modules we should be decomposing.  The modules $\wml$ 
studied here were defined by assigning a canonical choice of Drinfeld 
polynomials to any weight of the form $m\w_\l$.  It is unclear 
how this assignment should be generalized to all weights, or even if 
there is a correct generalization.

Chari and Pressley have investigated the notion of a ``minimal 
affinization'' of a $U_q(\g)$-module at some length, and report that 
while there is indeed a canonical representation of $U_q(\ghat)$ in 
our cases, often there is not.  For a generic highest weight of $D_n$ 
and $E_n$, in fact, there are three!  While the modules $\wml$ 
coincided with these ``minimal affinizations'' so far, there is no 
guarantee that the correct generalization would continue to do so.

\item
Finally, the polynomial relations among the characters should also be
generalized to highest weights other than multiples of a fundamental
weight.  Again, this would be interesting even in the $A_n$ case.
Preliminary investigation in this direction supports the guess that there
is a generalization of the form
$$
V(\lambda+\omega_\l)^2 =
V(\lambda) V(\lambda+2\omega_\l) +
\sum_k V(\mu_{i_k}) V(\mu_{j_k})
$$
where the weights $\mu$ in the sum of products are not very far away from
$\lambda$ and should have height at most that of $\lambda+\omega_\l$.
This again has the feeling of a discrete dynamical system about it.  If
one were to generalize the results of Chapter~\ref{chap_poly} to such a
system of relations, the set of ``base cases'' would grow at least to
$2^{\rank(\g)}-1$, since the recurrence relation would only give you
values for highest weights of the form $\lambda+2\omega_\l$.  Again, the
fact that there is not a canonical representation of $U_q(\g)$ with a
given highest weight makes things trickier.

\end{enumerate}

\end{document}